\documentclass{article}

\usepackage[nonatbib, main, final]{neurips_2026}

\usepackage[utf8]{inputenc} %
\usepackage[T1]{fontenc}    %
\usepackage{hyperref}       %
\usepackage{url}            %
\usepackage{booktabs}       %
\usepackage{amsfonts}       %
\usepackage{nicefrac}       %
\usepackage{microtype}      %

\usepackage{amsmath,amsthm,epsfig,amssymb,amsbsy,mathrsfs,mathtools}
\usepackage{enumerate}
\usepackage{bbm}
\usepackage{dsfont}
\usepackage[x11names]{xcolor}
\usepackage{algorithm}
\usepackage{algorithmic}
\usepackage{adjustbox}
\usepackage{subcaption}

\usepackage{wrapfig}

\usepackage{multirow,makecell}

\usepackage[capitalize,noabbrev,nameinlink]{cleveref}
\usepackage{crossreftools}
\usepackage[giveninits=true,maxbibnames=9, maxcitenames=2, backend=biber]{biblatex}
\usepackage{caption}
\usepackage{subcaption}
\usepackage{enumitem}
\usepackage{etoolbox}

\DeclareMathOperator*{\argmin}{arg\,min}

\DeclareMathOperator{\gph}{gph}

\DeclareMathOperator{\dom}{dom}

\DeclareMathOperator{\intr}{int}
\DeclareMathOperator{\bdry}{bdry}

\DeclareMathOperator{\id}{id}
\DeclareMathOperator{\arsinh}{arsinh}

\DeclareMathOperator{\zer}{zer}

\newcommand{\dOne}{\widetilde Fx^k}
\newcommand{\dTwo}{\widetilde Fz^k}

\newcommand{\pre}{\mathcal{K}}

\newcommand{\bR}{\mathbb{R}}

\newcommand{\bE}{\mathbb{E}}

\newcommand{\bN}{\mathbb{N}}

\newcommand{\exR}{\overline{\mathbb{R}}}

\newcommand{\cF}{\mathcal{F}}
\newcommand{\cC}{\mathcal{C}}

\newcommand{\cN}{\mathcal{N}}

\newcommand{\cG}{\mathcal{G}}

\makeatletter
\newcommand{\prox}[3][\@nil]{%
  \def\tmp{#1}%
   \ifx\tmp\@nnil
       \operatorname{prox}_{#3}^{#2}
    \else
         \operatorname{prox}_{#3}^{#1 \star #2}
    \fi}

\newcommand{\bprox}[3][\@nil]{%
  \def\tmp{#1}%
   \ifx\tmp\@nnil
       \operatorname{bprox}_{#3}^{#2}
    \else
        \operatorname{bprox}_{#3}^{#1 #2}
    \fi}
\makeatother

\usepackage{hyperref}
\hypersetup{
    colorlinks=true,
    linkcolor=Blue4,
    filecolor=magenta,
    citecolor = Blue4,  
    urlcolor=magenta
    }
\usepackage[capitalize]{cleveref}[0.19]

\crefname{section}{Section}{Sections}
\crefname{subsection}{Subsection}{Subsections}
\Crefname{section}{Section}{Sections}
\Crefname{subsection}{Subsection}{Subsections}

\Crefname{figure}{Figure}{Figures}

\crefformat{equation}{\textup{#2(#1)#3}}
\crefrangeformat{equation}{\textup{#3(#1)#4--#5(#2)#6}}
\crefmultiformat{equation}{\textup{#2(#1)#3}}{ and \textup{#2(#1)#3}}
{, \textup{#2(#1)#3}}{, and \textup{#2(#1)#3}}
\crefrangemultiformat{equation}{\textup{#3(#1)#4--#5(#2)#6}}%
{ and \textup{#3(#1)#4--#5(#2)#6}}{, \textup{#3(#1)#4--#5(#2)#6}}{, and \textup{#3(#1)#4--#5(#2)#6}}

\Crefformat{equation}{#2Equation~\textup{(#1)}#3}
\Crefrangeformat{equation}{Equations~\textup{#3(#1)#4--#5(#2)#6}}
\Crefmultiformat{equation}{Equations~\textup{#2(#1)#3}}{ and \textup{#2(#1)#3}}
{, \textup{#2(#1)#3}}{, and \textup{#2(#1)#3}}
\Crefrangemultiformat{equation}{Equations~\textup{#3(#1)#4--#5(#2)#6}}%
{ and \textup{#3(#1)#4--#5(#2)#6}}{, \textup{#3(#1)#4--#5(#2)#6}}{, and \textup{#3(#1)#4--#5(#2)#6}}

\newtheorem{theorem}{Theorem}[section]
\newtheorem{lemma}[theorem]{Lemma}
\newtheorem{definition}[theorem]{Definition}
\newtheorem{assumption}[theorem]{Assumption}
\crefname{assumption}{Assumption}{Assumptions}
\newtheorem{proposition}[theorem]{Proposition}
\newtheorem{corollary}[theorem]{Corollary}

\newtheoremstyle{boldremark}
    {\dimexpr\topsep/2\relax} 
    {\dimexpr\topsep/2\relax} 
    {}          
    {}          
    {\bfseries} 
    {.}         
    {.5em}      
    {}          

\theoremstyle{boldremark}
\newtheorem{remark}[theorem]{Remark}
\theoremstyle{boldremark}
\newtheorem{example}[theorem]{Example}

\newlist{thmenum}{enumerate}{1} 
\setlist[thmenum]{label=(\roman*), ref=\theproposition(\roman*), font=\rm}
\crefalias{thmenumi}{theorem} 

\AtEndEnvironment{remark}{\qed}
\AtEndEnvironment{example}{\qed}

\usepackage{xcolor}

\newcommand{\isocoe}{\kappa}

\title{Dual-Space Preconditioning for Variational Inequalities and Root-Finding Problems}

\author{Jan Quan\thanks{Equal contribution \\ \emph{E-mails:} 
 {\tt%
 \href{mailto:jan.quan@kuleuven.be}{\{jan.quan,} %
 \href{mailto:konstantinos.oikonomidis@kuleuven.be}{konstantinos.oikonomidis,} %
 \href{mailto:alexander.bodard@kuleuven.be}{alexander.bodard,} %
 \href{mailto:panos.patrinos@kuleuven.be}{panos.patrinos\}} %
\href{mailto:jan.quan@kuleuven.be,konstantinos.oikonomidis@kuleuven.be,alexander.bodard@kuleuven.be,panos.patrinos@kuleuven.be}{\!\!@kuleuven.be}}}
 \,,\;\, Konstantinos Oikonomidis$^*$,\;\, Alexander Bodard$^*$,\;\, Panagiotis Patrinos \\
KU Leuven (STADIUS)}

\begin{document}

\maketitle

\begin{abstract}
    This paper develops a dual-space preconditioning framework for variational inequalities and root-finding problems beyond classical cocoercivity and Lipschitz continuity assumptions. To this end, we introduce both a novel cocoercivity-type condition and a relaxed Lipschitz continuity condition, and derive convergence of deterministic and stochastic methods under these generalized conditions. For isotropic preconditioners, the resulting algorithms can be interpreted as variable-stepsize variants of classical methods, extending recent analyses based on $(L_0, L_1)$-type growth conditions while allowing simpler proofs and larger stepsize ranges. We also show that the framework naturally handles  variational inequalities through weak Minty-type arguments. Finally, we study the associated continuous-time dynamics and derive a variational and optimal control interpretation of the preconditioned flow. Overall, the results show that dual-space preconditioning provides a flexible mechanism for designing and analyzing operator methods beyond the standard Euclidean Lipschitz setting.
\end{abstract}

\section{Introduction}

In this paper, we study \emph{variational inequality} (VI) problems of the following form: 
\begin{equation*} \label{eq:VI} \tag{VIP}
    \text{find } x^\star \in C \text{ such that } \langle F x^\star, x - x^\star \rangle \geq 0
    \qquad \forall x \in C,
\end{equation*}
where $C$ is a nonempty, closed, convex subset of a Euclidean space $E$ with inner product $\langle \cdot,\cdot \rangle$ and induced norm $\|\cdot\|$, and $F : E \to E$ is a continuous operator. Throughout the paper, we assume that a solution $x^\star$ exists. Equivalently, \eqref{eq:VI} can be written as the inclusion problem
\(
    0 \in (F + N_C)(x^\star),
\)
where $N_C$ denotes the normal cone of $C$. %

Variational inequalities have received widespread attention \cite{burachik2008enlargements,facchinei2003finite,phelps1989convex} because they encompass problems such as (un)constrained minimization, multiplayer games, and min-max optimization. They arise in numerous machine learning applications, including reinforcement learning, generative adversarial network training, distributionally robust optimization, and optimal transport \cite{goodfellow2020generative,namkoong2016stochastic}. For instance, for the min-max problem
\(
    \min_{u\in U} \max_{v\in V} f(u,v),
\) where $U$ and $V$ are nonempty, closed, convex sets,
one can select the operator
\begin{equation} \label{eq:minimax-operator}
    F(u, v)
    =
    \bigl(\nabla_{u} f(u,v)^\top,
    -\nabla_{v} f(u,v)^\top \bigr)^\top .
\end{equation}
When $C=U\times V$ and $f$ is convex-concave, a solution to \eqref{eq:VI} then corresponds to a saddle point.

In recent years, problems of the form \eqref{eq:VI} have been studied under various monotonicity assumptions on $F$, ranging from strong monotonicity to Minty variational inequalities \cite{loizou2021stochastic,gorbunov2022stochastic,diakonikolas2021efficient,choudhury2023single,pethick2023escaping}. However, with the exception of some recent works on $\alpha$-symmetric $(L_0,L_1)$-Lipschitz continuity \cite{vankov2024generalized,vankov2024generalized2,choudhury2025extragradient}, most existing analyses rely on standard Lipschitz continuity of $F$. This contrasts with a growing trend in optimization toward generalized smoothness notions, which provide a more accurate description of practical problems \cite{zhang2019gradient,ahn2023linear} and lead to sharper algorithmic guarantees \cite{gorbunov2024methods}.

A related line of work concerns dual-space, or nonlinear, preconditioning methods \cite{maddison2021dual,laude2019optimization}. These methods are naturally analyzed under generalized smoothness assumptions \cite{laude2025anisotropic} and have proved useful in optimization, encompassing, for instance, gradient-clipping-like methods \cite{oikonomidis2025nonlinearlypreconditionedgradientmethods}. However, to date this perspective has largely remained restricted to optimization problems, leaving its potential for variational inequalities and root-finding problems underexplored. This motivates our central question:

\medskip
\emph{Can preconditioning techniques be extended to variational inequalities and root-finding problems, and can the resulting algorithms be analyzed under generalized Lipschitz continuity conditions?}
\medskip

We highlight two key challenges in this regard. First, classical cocoercivity and Lipschitz continuity assumptions must be generalized in a way that is compatible with the nonlinear preconditioner. However, existing results rely heavily on $F$ being the gradient of a function that plays an explicit role in the corresponding notion of anisotropic smoothness \cite{laude2025anisotropic}. Second, there is an inherent trade-off between the generality of the preconditioning framework and the strength of the resulting convergence guarantees. It is therefore important to identify meaningful classes of preconditioners that admit a unified analysis, cover a broad range of algorithms, and recover known results as special cases.

\paragraph{Contributions}
Our contributions can be summarized as follows.
\begin{itemize}[leftmargin=15pt]
    \item We introduce a novel cocoercivity-type condition (\cref{def:aniso_cocoer}) that generalizes classical cocoercivity through a preconditioner $\nabla \phi^*$. We relate this condition to dual relative smoothness~\cite{maddison2021dual} and use it to establish convergence results for preconditioned forward methods, and obtain last-iterate convergence rates for a broad class of isotropic preconditioners.

    \item We propose a generalized Lipschitz continuity condition (\cref{def:aniso_cont}) that encompasses and extends recent notions such as $\alpha$-symmetric $(L_0,L_1)$-Lipschitz continuity \cite[Assumption 2.1]{vankov2024generalized}. We use this condition to analyze preconditioned extragradient methods and, by focusing on isotropic preconditioners, derive variable-stepsize variants of several existing algorithms, including Extragradient+ (EG+) \cite{diakonikolas2021efficient} and Optimistic Gradient Descent Ascent (OGDA) \cite{kotsalis2022simple}. We further extend the analysis to stochastic settings with large batch sizes. Overall, our unified preconditioning framework offers a clear route to generalizing existing methods while substantially simplifying their convergence analyses.

    \item Finally, we study the associated continuous-time dynamics, develop a variational perspective on nonlinear preconditioning through the Fitzpatrick function \cite[p.~761]{borwein2016maximal}, and establish a fundamental connection with optimal control problems.
\end{itemize}

\paragraph{Notation}
We denote $\exR = \bR \cup\{\pm\infty\}$ and $\bN = \{0, 1, 2, \ldots\}$. For a proper convex function $f:E\to\exR$, its \emph{subdifferential} at $x\in E$ is $\partial f(x) = \{v \in E \mid f(\bar x) \geq f(x) + \langle v, \bar x - x\rangle, \forall \bar x\in E\}$. The \emph{indicator function} of a set $C \subseteq E$ is defined as $\delta_C(x) = 0$ if $x \in C$ and $+\infty$ otherwise. For a nonempty, closed, convex set $C$, its \emph{normal cone} is $N_C = \partial \delta_C$, that is, $N_C(\bar x) = \{y \in E \mid \langle y,x-\bar x \rangle \leq 0, \forall x \in C\}$ for $\bar x \in C$ and $N_C(\bar x)=\emptyset$ otherwise. Further, we denote by $P_C$  the \emph{Euclidean projection} onto $C$. We say that a proper, lsc, convex function $f:E \to \exR$ is \emph{essentially smooth} if $\intr (\dom f)\neq \emptyset$ and $f$ is differentiable on $\intr (\dom f)$ such that $\|\nabla f(x^\nu)\| \to \infty$, whenever $\intr (\dom f) \ni x^\nu \to x \in \bdry \dom f$, and \emph{essentially strictly convex}, if $f$ is strictly convex on every convex subset of $\dom \partial f \coloneq \{x \in E \mid \partial f(x) \neq \emptyset\}$, and \emph{Legendre}, if $f$ is both essentially smooth and essentially strictly convex. In particular, a smooth convex function on $E$ is essentially smooth. For a convex function $f \in \cC^1(E)$, we denote $D_f(x, \bar x) = f(x) - f(\bar x) - \langle \nabla f(\bar x), x - \bar x\rangle$ for $x,\bar x\in E$. For a set-valued mapping $T: E \rightrightarrows E$, we define its \emph{graph} as $\gph T \coloneq \{(x,y) \in E \times E \mid y \in T(x)\}$ and its \emph{zeros} $\zer T = \{x \in E \mid 0\in Tx\}$. We denote by $\id$  the \emph{identity operator} on the relevant space and write $\|\cdot\|_{\rm op}$ for the \emph{operator norm}.

\subsection{Related work}

\paragraph{Variational inequalities} In recent years, the VI literature has mostly been centered around lifting the monotonicity assumption, mainly due to applications in (non)convex-(non)concave minimax problems~\cite{lin2020gradient}. Particular attention has been given to a class of nonmonotone problems that satisfy the so-called \emph{Minty variational inequality} \cite{mertikopoulos2018optimistic,song2020optimistic,boct2021minibatch,kotsalis2022simple}. Its generalization, the \emph{weak Minty variational inequality}, was proposed in \cite{diakonikolas2021efficient} along with a convergent extragradient scheme, while \cite{pethick2023escaping} introduced an adaptive extension with increased stepsize ranges. Moreover, a large body of work has been devoted to stochastic extensions of classical methods \cite{hsieh2020explore,bohm2022two,gorbunov2023high,pethick2023solving,gorbunov2022clipped}.

\paragraph{Generalized smoothness}
Our work is closely related to notions of generalized smoothness that have received widespread attention recently in the functional setting. The most well-known condition, called $(L_0, L_1)$-smoothness~\cite{zhang2019gradient}, which for a twice differentiable objective $f$ takes the form $\|\nabla^2 f(x)\|_{\rm op} \leq L_0+L_1\|\nabla f(x)\|$ \cite{zhang2019gradient,hubler2024parameter}, allows the norm of the Hessian to grow unbounded. Various existing
methods have been analyzed under this new smoothness
condition \cite{wang2023convergence,faw2023beyond,koloskova2023revisiting}, and new methods have also been proposed \cite{gorbunov2024methods,vankov2024optimizing}.

\paragraph{Dual-space preconditioning}
In the \emph{functional} setting, the nonlinear preconditioning framework that we build upon was introduced in \cite{maddison2021dual} for gradient-type methods and in \cite{laude2019optimization} for proximal methods. In \cite{laude2025anisotropic}, the anisotropic proximal gradient algorithm for nonconvex composite minimization problems was introduced and studied under a new condition called anisotropic smoothness \cite{laude2023dualities}. Subsequent work showed that this condition encompasses $(L_0,L_1)$-smoothness and connected the framework to algorithms for neural network training \cite{oikonomidis2025nonlinearlypreconditionedgradientmethods,oikonomidis2025nonlinearly,bodard2025escaping,oikonomidis2026constrained}. The dual-space preconditioning method was also linked to the Lion algorithm \cite{chen2023symbolic} and to the solution of a certain constrained minimization problem in \cite{chen2023lion}, under Lipschitz smoothness of the objective. In the \emph{operator} setting, \cite{laude2023anisotropic} studied the convergence of a nonlinearly preconditioned \emph{backward} scheme for monotone operators, obtaining convergence guarantees under general assumptions. For a well-behaved family of preconditioners, a nonlinear proximal point scheme was already analyzed in \cite{luque1987nonlinear}.

\section{Convergence under \texorpdfstring{$\nabla \phi^*$}{∇ϕ∗}-cocoercivity} \label{sec:problem_form}
Throughout this section, we focus on preconditioners that are gradients of convex functions. Although we work in finite-dimensional spaces, to remain aligned with the classical convention that $F$ takes values in the dual space, denoting $\phi:E \to \exR$ for the so-called \emph{reference function}, we use the gradient of its convex conjugate, $\nabla\phi^*$, as the preconditioner.
\begin{assumption} \label{assum:phistar}
    $\phi^*:E \to \bR$ is even and Legendre and satisfies $\phi^* \geq 0$ and $\phi^*(0) = 0$.
\end{assumption}

Consider the \emph{root-finding} problem, i.e., \eqref{eq:VI} with $C=E$,
\begin{equation*} \tag{RFP} \label{eq:rfp}
    \text{find } x^\star \in E \text{ such that } Fx^\star = 0.
\end{equation*}
In the standard Euclidean setup, the simplest assumption that leads to convergence of \emph{forward} schemes for \eqref{eq:rfp} is \emph{cocoercivity} \cite[Definition 4.10]{bauschke2017correction} of $F$. This condition covers a large class of monotone operators, including  gradients of convex and Lipschitz smooth functions. We propose the following generalization that is based on the notion of $\nabla \phi^*$-firm nonexpansiveness from \cite[Definition 4.1]{wang2022bregman}.
\begin{definition}[$\nabla\phi^*$-cocoercivity] \label{def:aniso_cocoer}
    Let $\beta > 0$. We say that $F$ is $\beta$-$\nabla \phi^*$-cocoercive if the following inequality holds for all $x, \bar x \in E$:
    \begin{equation} \label{eq:aniso_coco}
        \beta \langle Fx-F\bar x, \nabla \phi^*(Fx) - \nabla \phi^*(F \bar x)\rangle \leq \langle F x - F\bar x, x - \bar x \rangle.
    \end{equation}
\end{definition}

\begin{figure}
    \centering
    \begin{subfigure}[b]{0.32\textwidth}
        \includegraphics[width=\linewidth]{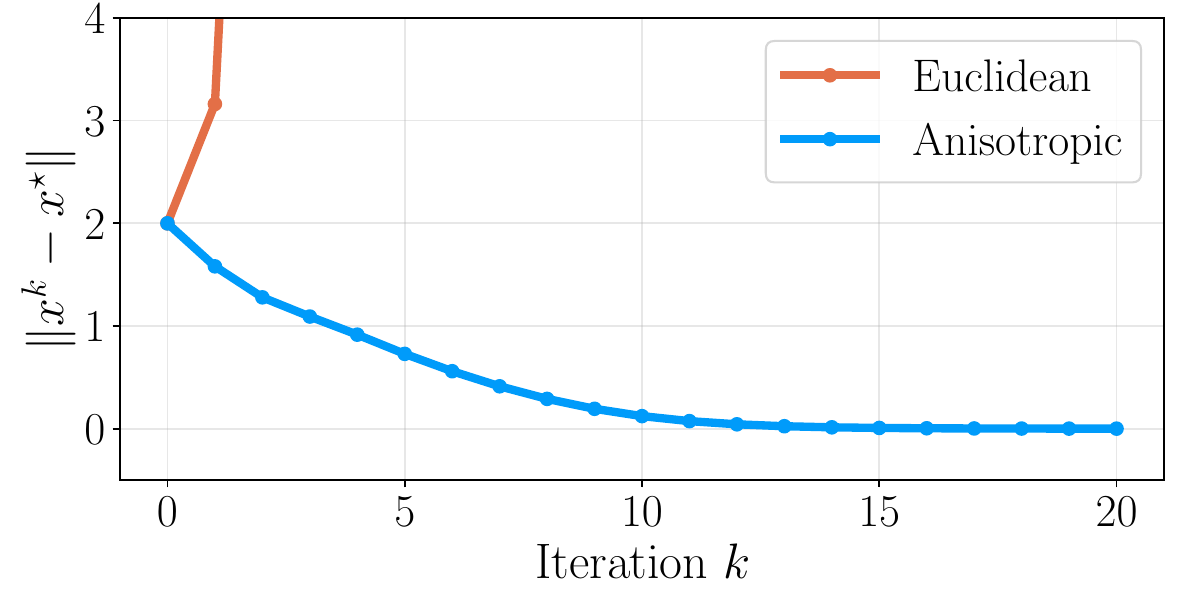}
    \end{subfigure} \hfill
    \begin{subfigure}[b]{0.32\textwidth}
        \includegraphics[width=\linewidth]{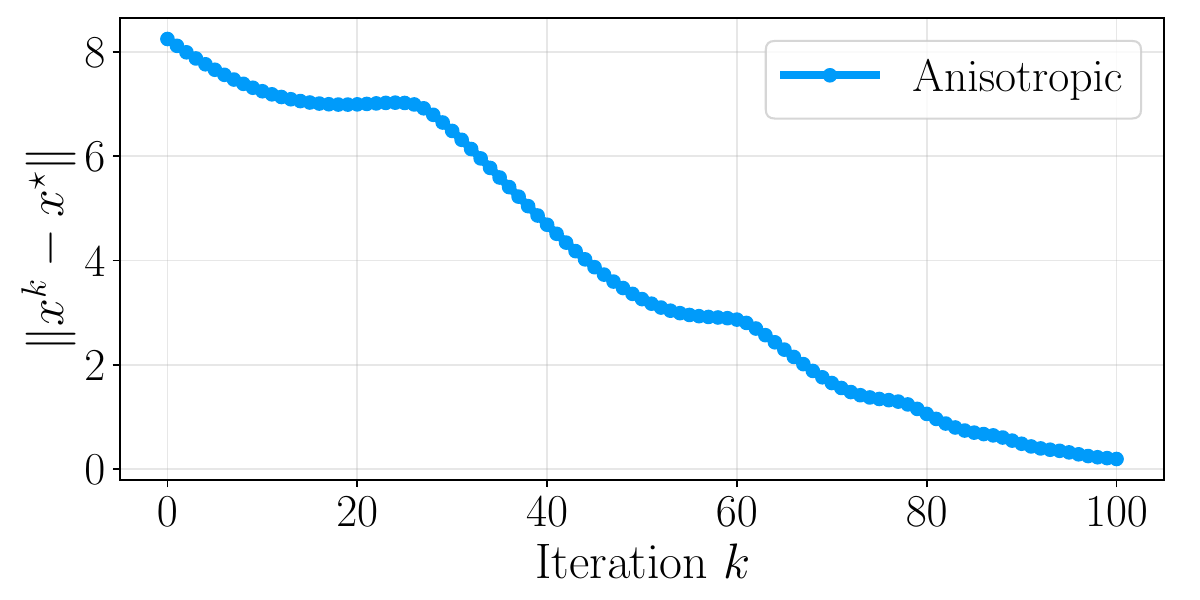}
    \end{subfigure} \hfill
    \begin{subfigure}[b]{0.32\textwidth}
        \includegraphics[width=\linewidth]{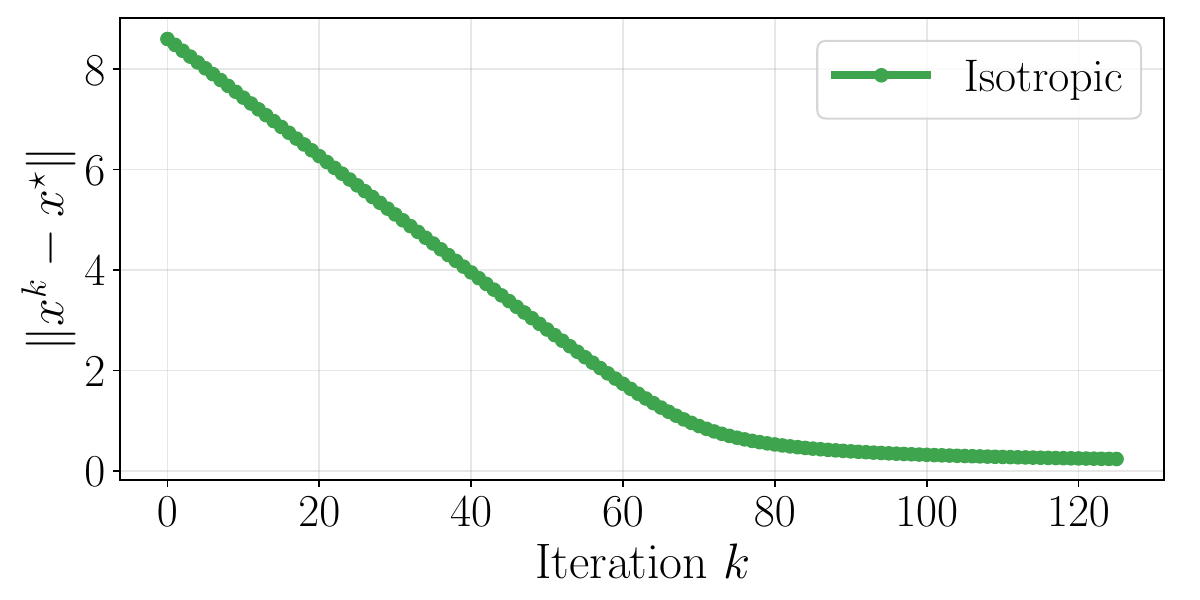}
    \end{subfigure} 
    \caption{(left) Comparison of the Euclidean forward step and \eqref{eq:GD} with $\gamma = \tfrac{1}{2}$, and $\phi^*$ and $F$ as in \cref{ex:aniso_cocoer}. (middle) Convergence of \eqref{eq:GD} with $\gamma = \tfrac{1}{5}$, and $\phi^*$ and $F$ as in \cref{ex:aniso_not_fejer}. The sequence is not necessarily Fej\'er monotone for non-isotropic preconditioners. (right) Convergence of \eqref{eq:GD} with $\gamma = \frac16$, and $\phi^*$ and $F$ as in \cref{ex:iso_coco}.}
    \label{fig:toy_examples}
\end{figure}

For $\phi^* = \tfrac{1}{2}\|\cdot\|^2$, this definition reduces to classical cocoercivity. More generally, $\nabla\phi^*$-cocoercivity extends the dual relative smoothness condition of \cite[Assumption 3.1.3]{maddison2021dual} to operators, as the following proposition shows.
\begin{proposition}[$\nabla\phi^*$-cocoercivity from dual relative smoothness] \label{thm:dual_relative}
    Suppose that $\phi^*$ satisfies \cref{assum:phistar}.
    Let $f \in \cC^1(E)$ be convex and such that the following inequality holds for all $x, \bar x \in E$:
    \begin{equation} \label{eq:smooth_relative}
        D_{\phi^*}(\nabla f(x), \nabla f(\bar x)) \leq L D_f(\bar x, x), \qquad L > 0.
    \end{equation}
    Then, $\nabla f$ is $L^{-1}$-$\nabla \phi^*$-cocoercive.
\end{proposition}
The natural way to tackle \eqref{eq:rfp} under \cref{def:aniso_cocoer} is then through the following generalization of the dual-space preconditioning method from \cite{maddison2021dual}:
\begin{equation} \label{eq:GD}
    x^{k+1} = x^k - \gamma \nabla \phi^*(Fx^k).
\end{equation}
The analysis of \eqref{eq:GD} when $\phi$ has full domain can in fact be obtained by rewriting the scheme as the nonlinear proximal point iteration studied in \cite{laude2023anisotropic}: $x^{k+1} = (\id + \gamma\nabla \phi^* \circ T)^{-1}(x^k)$ for some maximally monotone operator $T$. To make the paper self-contained, we provide a direct proof based on the forward step in \cref{subsec:gd_res}.
\begin{theorem} \label{thm:gd_general_results}
    Suppose that $\phi^*$ satisfies \cref{assum:phistar}.
    Let $F$ be $\nabla \phi^*$-cocoercive with constant $\beta$. Let, moreover, $\{x^k\}_{k \in \bN}$ be the sequence of iterates generated by \eqref{eq:GD} with $\gamma\in(0,\beta]$. Then, $\phi^*(Fx^{k+1}) \leq \phi^*(Fx^k)$ for all $k \in \bN$ and $\{D_{\phi^*}(Fx^{k}, Fx^{k+1})\}_{k=0}^\infty$ is summable. In particular, $\lim_{k \to \infty} D_{\phi^*}(Fx^{k}, Fx^{k+1}) = 0$. Moreover, $\{Fx^{k}\}_{k \in \bN}$ is bounded and every limit point of $\{x^k\}_{k \in \bN}$ is a zero of $F$.
\end{theorem}
\begin{example} \label{ex:aniso_cocoer}
    Consider the dual reference function $\phi^*(y_1, y_2) = h^*(y_1) + h^*(y_2)$ with $h^*(y) = |y| - \log(|y| + 1)$ and the operator $F(x_1, x_2) = \begin{psmallmatrix}
        x_1 + x_1^3 \\ x_1 + x_2
    \end{psmallmatrix}$. 
    This $F$ is not standard Euclidean cocoercive, in fact not even Lipschitz, while it is $\nabla\phi^*$-cocoercive with constant $\beta=\tfrac{1}{2}$. Moreover, for any stepsize $\gamma > 0$, there exists some initial point so that the Euclidean forward step diverges. In contrast, the initial sublevel set of $\phi^* \circ F$ is bounded, so convergence of \eqref{eq:GD} is guaranteed for $0 < \gamma \leq \beta$ by the preceding result. The left subfigure of \Cref{fig:toy_examples} shows divergence of the Euclidean forward step, while \eqref{eq:GD} converges.
\end{example}
\begin{example} \label{ex:aniso_not_fejer}
    Consider the dual reference function $\phi^*(y_1, y_2) = h^*(y_1) + h^*(y_2)$ with $h^*(y) = |y| - \log(|y| + 1)$ and the operator $F(x_1, x_2) = \begin{psmallmatrix}
        x_1 - 2 x_2 \\ 2x_1 + x_2
    \end{psmallmatrix}$. 
    It can be shown that $F$ is $\nabla\phi^*$-cocoercive with constant $\beta=\tfrac{1}{5}$. In contrast to \cref{ex:aniso_cocoer}, the iteration is not Fej\'er monotone, as shown in the middle subfigure of \Cref{fig:toy_examples}.
\end{example}

When $\nabla \phi^*(y) = \isocoe(\|y\|)y$ for some $\kappa : [0, +\infty) \to (0, +\infty)$, we say that $\nabla\phi^*$ is \emph{isotropic}. Isotropy yields stronger convergence guarantees and permits larger stepsizes, as shown by the following theorem, which also provides last-iterate convergence rates.
\begin{theorem}[Convergence for isotropic preconditioners] \label{thm:iso_conv_GD}
    Let $F$ be $\beta$-$\nabla \phi^*$-cocoercive. Let, moreover, $\{x^k\}_{k \in \bN}$ be generated by \eqref{eq:GD} with an isotropic preconditioner $\nabla \phi^*$, where $\phi^* = h^*\circ\|\cdot\|$ and $h^*$ satisfies \cref{assum:phistar}. Then, the following hold  for every $x^\star\in\zer F$ and $K\in\bN$:
    \begin{thmenum}[leftmargin=*, labelindent=6pt]
        \item \label{thm:cocoer:primal_space_fejer} For $\gamma \in (0, 2\beta)$, $\{x^k\}_{k \in \bN}$ is bounded and Fej\'er monotone with respect to $\zer F$, that is, $\|x^{k+1} - x^\star\| \leq \|x^k - x^\star\|$ for all $x^\star\in \zer F$ and $k\in\bN$.
        \item \label{thm:cocoer:rate} For $\gamma \in (0,\beta]$, $\{\|Fx^k\|\}_{k\in \bN}$ is a nonincreasing sequence and we have the following rate:
         \begin{equation*}
            \|\nabla \phi^*(Fx^{K})\|^2 \leq \frac{\|x^0-x^\star\|^2}{(K+1) \gamma(2\beta - \gamma)}.
        \end{equation*}
        If, moreover, $\isocoe$ is a nonincreasing function on $[0,+\infty)$, we obtain the last-iterate rate
        \[
            \|Fx^{K}\| \leq \frac{\|x^0-x^\star\|}{\isocoe(\|Fx^0\|)\sqrt{\gamma(2\beta - \gamma)(K+1)}}.
        \]
    \end{thmenum}
\end{theorem}
\begin{example} \label{ex:iso_coco}
   Consider the dual reference function $\phi^*(y) = \|y\| - \log(\|y\| + 1)$ and the operator $F(x) = \|x\|^2(x + Jx)$ with $J$ a linear operator that satisfies $J^*J = \id$ and $J^* = -J$, where $J^*$ denotes the adjoint of $J$, e.g., $J = \begin{psmallmatrix}
        0 & -1 \\ 1 & 0
    \end{psmallmatrix}$ in $2$ dimensions. It can be shown that $F$ is $\nabla\phi^*$-cocoercive with constant $\beta=\frac{1}{6}$ but is not cocoercive in the standard Euclidean sense. The right subfigure of \Cref{fig:toy_examples} shows convergence of \eqref{eq:GD}. Note that the iterates are Fej\'er monotone, as described in \cref{thm:iso_conv_GD}.
\end{example}

Another example in a similar spirit to \cref{ex:iso_coco} for port-Hamiltonian equilibrium operators is provided in \cref{app:port_hamiltonian}.

\subsection{Extension to the composite setting} \label{subsec:extension}
We now consider the more general inclusion problem
\begin{equation*}
    0 \in Fx + Gx,
\end{equation*}
where $G:E \rightrightarrows E$ is maximally monotone, subsuming \eqref{eq:VI} by choosing $G=N_C$. Moreover, we assume that the reference function $\phi$ has full domain. The following forward-backward algorithm extends the method of \cite{laude2025anisotropic} to the operator setting:
\begin{equation} \label{eq:for_back}
    y^k = x^k - \gamma \nabla \phi^*(Fx^k), \qquad x^{k+1} = (\id + \gamma \nabla \phi^* \circ G)^{-1}(y^k).
\end{equation}
In light of \cite[Proposition 3.5]{laude2023anisotropic}, the operator $(\id + \gamma \nabla \phi^* \circ G)^{-1}$ is single-valued. Our aim is to write \eqref{eq:for_back} as a single forward step, so that we can directly use \cref{thm:gd_general_results}. Define the forward-backward map $T_\gamma(x)\coloneq (\id + \gamma \nabla \phi^* \circ G)^{-1}(x-\gamma \nabla \phi^*(Fx))$ and the transformed fixed-point residual $R_\gamma(x) \coloneq \nabla \phi(\tfrac{x-T_\gamma(x)}{\gamma})$. Note that $\zer R_\gamma = \zer (F + G)$ and that \eqref{eq:for_back} can be written as $x^{k+1} = x^k - \gamma \nabla \phi^*(R_\gamma(x^k))$. We are thus led to the following corollary.
\begin{corollary}
    Suppose that $\phi^*$ satisfies \cref{assum:phistar}, $\dom\phi=E$, and $G:E\rightrightarrows E$ is maximally monotone. Let $\gamma > 0$ and suppose that $R_\gamma$ is $\gamma$-$\nabla\phi^*$-cocoercive. Let, moreover, $\{x^k\}_{k\in\bN}$ be the sequence of iterates generated by \eqref{eq:for_back}. Then, $\phi^*(R_\gamma (x^{k+1})) \leq \phi^*(R_\gamma (x^k))$ for all $k \in \bN$, and $\{D_{\phi^*}(R_\gamma (x^k), R_\gamma (x^{k+1}))\}_{k=0}^\infty$ is summable. In particular, $\lim_{k \to \infty} D_{\phi^*}(R_\gamma (x^k), R_\gamma (x^{k+1})) = 0$. Moreover, $\{R_\gamma (x^k)\}_{k \in \bN}$ is bounded and every limit point of $\{x^k\}_{k \in \bN}$ is a zero of $F + G$.
\end{corollary}
When $G\equiv 0$, we have $R_\gamma = F$, and we recover \eqref{eq:GD} and \cref{thm:gd_general_results}. When $F\equiv 0$, we recover the method of \cite{laude2023anisotropic}, which we discuss further in \cref{app:discuss_composite}. In the truly composite case, our analysis requires $\nabla\phi^*$-cocoercivity of $R_\gamma$, and we currently do not have separate conditions on $F$ and $G$ that imply this property. Moreover, several natural Lyapunov candidates need not decrease along the iterates when $F$ is $\nabla\phi^*$-cocoercive and $G$ is monotone. In particular, we numerically tested $D_{\phi^*}(v^{k+1}, -Fx^k)$, $D_{\phi^*}(-Fx^k,v^{k+1})$, $\tfrac{1}{2}\|x^k-x^\star\|^2$, $D_{\phi^*}(Fx^k,Fx^\star)+D_{\phi^*}(v^{k+1},v^\star)$, and $\phi^*(\nabla \phi(\nabla \phi^*(Fx^k)+\nabla \phi^*(v^{k+1})))$, where $v^{k+1} \in Gx^{k+1}$ and $v^\star = - Fx^\star$ with $x^\star \in \zer(F + G)$. We found a simple example in which each quantity increases during at least one of the first ten iterations, although the iterates numerically approach a solution. This illustrates the difficulty of obtaining a general Lyapunov analysis for \eqref{eq:for_back}. The code used to generate the example is included in the supplementary material.

\section{Extragradient methods under preconditioned Lipschitz continuity} \label{sec:prec_lip}
In this section, we focus on isotropic preconditioners and as such denote them by $\pre$ instead of $\nabla\phi^*$.
\begin{assumption} \label{assum:isotropy}
    The preconditioner $\pre:E \to E$ is monotone and continuous, satisfies $\pre(0) = 0$, and is isotropic: $\pre(y) = \kappa(\|y\|)y$ for all $y\in E$ and some $\kappa:[0,+\infty) \to (0, +\infty)$.
\end{assumption}
Such a preconditioner can be identified with the gradient of the dual reference function $\phi^*(y) = \int_0^{\|y\|} t\kappa(t)\,dt$, for which $\nabla\phi^*=\pre$. Note that we have relaxed the strict monotonicity, which is implied by \cref{assum:phistar}, to merely monotonicity.

\subsection{The preconditioned Lipschitz continuity assumption}

Extragradient methods are usually analyzed under Lipschitz continuity of $F$. Instead, we impose Lipschitz continuity on the preconditioned operator $\pre\circ F$.
\begin{assumption}[$\pre$-Lipschitz continuity] \label{def:aniso_cont}
    Let $L > 0$. We say that $F : E \to E$ is \(\pre\)-Lipschitz continuous with constant $L$ if the following inequality holds for all $x, \bar x \in E$: 
    \begin{equation} \label{eq:aniso_cont}
        \| \pre(Fx) - \pre(F \bar x)\| \leq L\|x-\bar x\|.
    \end{equation}
\end{assumption}
This condition is a direct generalization of \cite[Assumption 2.5]{oikonomidis2025nonlinearly} to the operator setting.
If $\pre$ is Lipschitz continuous, every Lipschitz-continuous operator $F$ satisfies \cref{def:aniso_cont} trivially. In fact, \cref{def:aniso_cont} is significantly \emph{weaker} than classical Lipschitz continuity for judicious choices of $\pre$. To see this, suppose that $F \in \cC^1(E)$ satisfies the $(L_0, L_1)$-Lipschitz continuity bound \cite{zhang2019gradient} $\forall x \in E: \| \nabla F(x) \|_{\rm op} \leq L_0 + L_1 \| Fx \|$ with $L_0>0,L_1 > 0$. Then, it can be shown by arguments similar to those in \cite[Proposition 2.6]{oikonomidis2025nonlinearly} that for $\phi^*(y) = \frac{\|y\|}{L_1} - \frac{L_0}{L_1^2} \ln(1 + \frac{L_1}{L_0} \|y\|)$, $\forall x \in E: \| \nabla (\pre\circ F)(x) \|_{\rm op} \leq 1$, with $\pre = \nabla \phi^*$, which is a sufficient condition for \eqref{eq:aniso_cont} with $L=1$.
More generally, the class of $\alpha$-symmetric $(L_0, L_1)$-Lipschitz continuous operators satisfies \cref{def:aniso_cont} for a suitable preconditioner $\pre$, as we show in \cref{prop:alpha_sym_smooth}.
\begin{definition}[$\alpha$-symmetric $(L_0, L_1)$-Lipschitz continuity \cite{vankov2024generalized}] \label{def:alpha_sym_smooth}
    We say that $F:E\to E$ is $\alpha$-symmetric $(L_0, L_1)$-Lipschitz for some $L_0 > 0$, $L_1 \geq 0$ and $\alpha\in(0,1]$ if the following inequality holds for all $x,\bar x\in E$:
    \begin{equation} \label{eq:alpha_sym_smooth}
        \|Fx - F\bar x\| \leq (L_0 + L_1 \max_{\theta\in[0,1]} \|F(\theta x + (1-\theta)\bar x)\|^\alpha)\|x - \bar x\|.
    \end{equation}
\end{definition}
\begin{proposition} \label{prop:alpha_sym_smooth}
    If $F : E \to E$ satisfies \Cref{def:alpha_sym_smooth} for some $L_0 > 0$, $L_1\geq 0$ and $\alpha\in(0,1]$, then  $F$ is $\pre$-Lipschitz continuous with constant $L = 1$ and $\pre : y \mapsto \frac{y}{L_0 + L_1\|y\|^\alpha}$.
\end{proposition}
It is natural to ask whether the converse also holds. The following example shows that this is in fact not the case, indicating that \cref{def:aniso_cont} is more general.
\begin{example} \label{ex:diverging-toilet}
    Consider $F:\bR^2 \to \bR^2$ defined by $Fx = e^{\|x\|^2/2}\begin{psmallmatrix}
        -x_2\\ x_1
    \end{psmallmatrix}$. Then, $F$ is $\pre$-Lipschitz continuous with constant $1$ and $\pre : y \mapsto \frac{y}{1+\|y\|}$, while it is not $\alpha$-symmetric $(L_0,L_1)$-Lipschitz for any $L_0,L_1,\alpha$ satisfying the requirements of \Cref{def:alpha_sym_smooth}. Similar examples in the functional setup are given in \cite[\S A.3]{bodard2025escaping}.
\end{example}
A partial converse does hold, however, when $0 < \alpha < 1$ in \cref{prop:alpha_sym_smooth}.
\begin{remark} \label{rem:alpha_sym_conv}
    If $0 < \alpha < 1$ and $F$ is $\pre$-Lipschitz continuous with constant $L>0$ and $\pre$ as in \cref{prop:alpha_sym_smooth}, then $F$ is $\alpha$-symmetric $(LL_0/(1-\alpha), LL_1/(1-\alpha))$-Lipschitz continuous.
\end{remark}

\subsection{Convergence under preconditioned Lipschitz continuity} \label{sec:convergence-preconditioned-lipschitz}

This subsection describes in detail how existing convergence results for extragradient and related methods under classical Lipschitz continuity can be generalized to hold under $\pre$-Lipschitz continuity.
To that end, we first recall the weak Minty variational inequality.
\begin{assumption} \label{assump:generalized-wmvi}
    The operator $\pre \circ F + N_C$ satisfies the weak Minty variational inequality (MVI), i.e., there exists a nonempty set $S^\star \subseteq \zer (\pre \circ F + N_C)$ and some $\rho \in \bR$ such that for all $x^\star \in S^\star$
    \begin{equation} \label{eq:wmvi}
        \langle y,x-x^\star \rangle \geq \rho \|y\|^2 \qquad \forall (x,y) \in \gph (\pre \circ F + N_C).
    \end{equation}
\end{assumption}
The zeros of the operator $\pre\circ F+N_C$ appearing in \cref{assump:generalized-wmvi} are exactly those of $F + N_C$.
\begin{lemma} \label{lem:zero-correspondance}
    If $\pre$ satisfies \cref{assum:isotropy}, then $\zer (\pre\circ F + N_C) = \zer (F + N_C)$.
\end{lemma}
The following key result establishes that monotonicity of $F$ implies that the operator $\pre \circ F + N_C$ satisfies the weak MVI with $\rho = 0$.
\begin{proposition} \label{lem:T-wmvi}
    If $\pre$ satisfies \cref{assum:isotropy} and \(F : E \to E\) is monotone, then \cref{assump:generalized-wmvi} holds with $S^\star = \zer (F + N_C)$ and $\rho = 0$.
\end{proposition}
\begin{remark}
    Suppose that $C=E$ and that there exists a nonempty set $S^\star\subseteq\zer F$ and $\rho\in\bR$ such that, for every $x^\star\in S^\star$
    \begin{equation} \label{eq:aniso-wmvi}
        \langle Fx, x - x^\star \rangle \geq \rho \langle \pre(Fx), Fx\rangle, \qquad \forall x \in E.
    \end{equation}
    It is easily verified that \eqref{eq:aniso-wmvi} is equivalent to $\pre \circ F$ satisfying the weak MVI with constant $\rho$.
    If $\rho < 0$, then this condition may be stricter or less strict than the usual weak MVI, depending on $\kappa$. Notably, \eqref{eq:aniso-wmvi} is an extension of \eqref{eq:aniso_coco} centered at the point $\bar x = x^\star$ and for negative constants $\beta$, generalizing again the classical setup and the relation between the weak MVI and cohypomonotonicity \cite{evens2023convergence}.
\end{remark}

The combination of \cref{def:aniso_cont,assump:generalized-wmvi} turns out extremely powerful, as they allow for the immediate transfer of convergence results under standard Lipschitz continuity of $F$ to the preconditioned Lipschitz continuity setting by formally replacing $F$ with $\pre \circ F$.
We now demonstrate this transfer of convergence results for two extragradient-type methods.

First, consider a preconditioned version of the CEG+ algorithm from \cite{pethick2023escaping}, for simplicity written with fixed stepsizes. 
Given the forward operator $H \coloneq \id - \gamma \pre \circ F$, it performs iterations
\begin{equation*} \label{eq:EG+} \tag{P-CEG+}
    z^k = P_C(x^k - \gamma \pre(Fx^k)), \qquad x^{k+1} = x^k + \beta (H z^k-H x^k).
\end{equation*}
By selecting $\pre= \id$ we recover the original CEG+ algorithm \cite{pethick2023escaping}, and by further restricting to the unconstrained setting where $C = E$ this reduces to the EG+ method of \cite{diakonikolas2021efficient} since $Hz^k - H x^k = - \gamma \pre( F z^k ) = -\gamma Fz^k$.

Under \cref{def:aniso_cont,assump:generalized-wmvi}, $\pre \circ F$ is Lipschitz continuous and $\pre \circ F + N_C$ satisfies a weak MVI. Moreover, we have $\zer (\pre\circ F + N_C) = \zer (F + N_C)$ by \cref{lem:zero-correspondance}.
Since the normal cone $N_C$ is a maximally monotone operator, all hypotheses of \cite[Corollary 3.2]{pethick2023escaping} hold, directly yielding the following.

\begin{theorem} \label{th:p-ceg+-convergence}
    Suppose that \cref{assum:isotropy,def:aniso_cont,assump:generalized-wmvi} hold with $\rho > -\frac{1}{2L}$.
    Let $x^\star\in S^\star$, $\gamma \in ([ -2\rho ]_+, \frac{1}{L}]$ with $[-2\rho]_+ = \max\{-2\rho, 0\}$, and $\beta \in (0, 1 +
    \frac{2 \rho}{\gamma})$. Consider the sequences $\{x^k\}_{k \in \bN}$, $\{z^k\}_{k \in \bN}$ generated by \eqref{eq:EG+}. Then, for every $K\in\bN$,
    \begin{equation*}
        \min_{k=0,1,\dots,K}\|Hz^k-Hx^k\|^2 \leq \frac{\|x^0-x^\star\|^2}{\sigma(K+1)},
    \end{equation*}
    where $\sigma = \beta(1+\frac{2\rho}{\gamma}-\beta)$. Furthermore, $\{z^k\}_{k \in \bN}$ is bounded and its limit points belong to $\zer(F + N_C)$.
\end{theorem}

In a similar manner, we can analyze methods that require a single evaluation of $F$ per iteration, including the following preconditioned variant of the Optimistic Gradient Descent Ascent method:
\begin{equation} \label{eq:p-ogda} \tag{P-OGDA}
    x^{k+1} = P_C(x^k - \gamma(2\pre(Fx^k)-\pre(Fx^{k-1}))), \qquad k \geq 1, \qquad x^0 = x^1 \in C.
\end{equation}
If the operator $F$ is monotone, then $\langle \pre (Fx), x - x^\star\rangle \geq 0$ for all $x\in C$ and $x^\star\in\zer(F+N_C)$ due to \cref{lem:T-wmvi} and choosing the zero normal vector in \eqref{eq:wmvi}.
Thus, all hypotheses of \cite[Theorem 2.5]{kotsalis2022simple} are satisfied, yielding the following result.
\begin{theorem}
    Let $F : E \to E$ be monotone and 
    suppose that \cref{assum:isotropy,def:aniso_cont} hold. Let $x^\star\in\zer(F + N_C)$.
    Consider the sequence $\{x^{k}\}_{k\in\mathbb{N}}$ generated by \eqref{eq:p-ogda}
    with $\gamma = \frac{1}{3L}$ and $x^0 = x^1 \in C$. Fix $K\geq1$. Let $R \in \argmin_{1\leq k\leq K} (\|x^{k+1} - x^k\|^2 + \|x^k - x^{k-1}\|^2)$. Then,
    \[
        \min_{y \in N_C(x^{R+1})} \|y + \pre(Fx^{R+1})\| \leq 20\sqrt{\frac{3}{2}} \frac{L \|x^{0} - x^\star\|}{\sqrt{K}}.
    \]
\end{theorem}

\subsection{Convergence under \texorpdfstring{$\alpha$}{α}-symmetric \texorpdfstring{$(L_0, L_1)$}{(L0, L1)}-Lipschitz continuity}

We now compare our framework with other adaptive stepsize rules for extragradient methods for solving \eqref{eq:rfp}, though we can also accommodate constraints in light of \cref{th:p-ceg+-convergence}. We focus on a modification of the Extragradient+ algorithm \cite{diakonikolas2021efficient} with adaptive stepsizes:
\begin{equation*}\label{eq:EG+-variable} \tag{EG+}
    z^k = x^k - \gamma_k F x^k, \qquad x^{k+1} = x^k - \omega_k F z^k.
\end{equation*}
For example, selecting $C = E$ and using $\pre(y) = \isocoe(\Vert y \Vert) y$, the iterates of \eqref{eq:EG+} correspond to those of \eqref{eq:EG+-variable} with
\(
    \gamma_k = \gamma \isocoe(\Vert F x^k \Vert)
\) and 
\(
    \omega_k = \omega \isocoe(\Vert F z^k \Vert)
\)
where $\omega = \beta \gamma$.
Other adaptive stepsizes have been analyzed in \cite{vankov2024generalized,choudhury2025extragradient} under $\alpha$-symmetric $(L_0, L_1)$-Lipschitz continuity.

\begin{table*}[t]
\centering
\small
\setlength{\tabcolsep}{4pt}
\renewcommand{\arraystretch}{1}
\caption{Summary of stepsize choices for \eqref{eq:EG+-variable} under \(\alpha\)-symmetric \((L_0,L_1)\)-Lipschitz continuity.}
\label{tab:stepsizes-comparison}
\begin{adjustbox}{max width=\textwidth}
\begin{tabular}{@{}c c c c c c@{}}
\toprule
Assumption & Reference & $(L_0,L_1)$ & \(\alpha\) & \(\gamma_k\) & \(\omega_k\) \\
\midrule
\multirow{5}{*}{\makecell{Strongly\\Monotone\textsuperscript{(1)}}}
& \cite{mokhtari2020unified} & \makecell{\(L\)-Lipschitz\textsuperscript{(4)}} & -- & \(\tfrac{0.25}{L}\) & \(\gamma_k\) \\
& \cite{vankov2024generalized} & \makecell{\((L_0,L_1)\)} & \(1\) & \(\displaystyle \min\!\left\{\tfrac{1}{4\mu},\tfrac{1}{2\sqrt{2}eL_0},\tfrac{1}{2\sqrt{2}eL_1\|F(x^k)\|}\right\}\) & \(\gamma_k\) \\
& \cite{choudhury2025extragradient} & \makecell{\((L_0,L_1)\)} & \(1\) & \(\displaystyle \tfrac{0.21}{L_0+L_1\|Fx^k\|}\) & \(\gamma_k\) \\
& \cite{choudhury2025extragradient} & \makecell{\((L_0,L_1)\)} & \((0,1)\) & \(\displaystyle \tfrac{0.61}{2K_0+\left(2K_1+2^{1-\alpha}K_2^{1-\alpha}\right)\|Fx^k\|^\alpha}\)\textsuperscript{(2)} & \(\gamma_k\) \\
& \textbf{Ours} & \makecell{\((L_0,L_1)\)} & $(0, 1]$ & $\tfrac{\gamma}{L_0+L_1 \Vert Fx^k\Vert^\alpha}$\textsuperscript{(3)} & $\frac{\omega}{L_0+L_1 \Vert Fz^k \Vert^\alpha}$\textsuperscript{(3)} \\
\midrule
\multirow{4}{*}{\makecell{Monotone\textsuperscript{(1)}}}
& \cite{gorbunov2022extragradient} & \makecell{\(L\)-Lipschitz\textsuperscript{(4)}} & -- & \(\tfrac{1}{\sqrt{2}L}\) & \(\gamma_k\) \\
& \cite{choudhury2025extragradient} & \makecell{\((L_0,L_1)\)} & \(1\) & \(\displaystyle \tfrac{0.45}{L_0+L_1\|Fx^k\|}\) & \(\gamma_k\) \\
& \cite{choudhury2025extragradient} & \makecell{\((L_0,L_1)\)} & \((0,1)\) & \(\displaystyle \tfrac{1}{2\sqrt{2}K_0+\left(2\sqrt{2}K_1+2^{3(1-\alpha)/2}K_2^{1-\alpha}\right)\|Fx^k\|^\alpha}\)\textsuperscript{(2)} & \(\gamma_k\) \\
& \textbf{Ours} & \makecell{\((L_0,L_1)\)} & $(0, 1]$ & $\tfrac{\gamma}{L_0+L_1 \Vert Fx^k \Vert^\alpha}$\textsuperscript{(3)} & $\frac{\omega}{L_0+L_1 \Vert Fz^k \Vert^\alpha}$\textsuperscript{(3)} \\
\bottomrule
\end{tabular}
\end{adjustbox}

\vspace{0.4em}
\begin{minipage}{0.98\textwidth}
\footnotesize
\textsuperscript{(1)} Convergence measures: strongly monotone: \(\|x^K-x^\star\|^2\); monotone: \(\min_{0\le k\le K}\|Fx^k\|^2\). %
\textsuperscript{(2)} For \(K_0,K_1,K_2\), see \cite[Proposition 3.1]{choudhury2025extragradient}. If \(L_1=0\), then \(K_1=K_2=0\).
\textsuperscript{(3)} $0 < \omega < \gamma < 1$.
\textsuperscript{(4)} Corresponds to $L_0 = L$ and $L_1 = 0$.
\end{minipage}
\end{table*}

\begin{wrapfigure}{r}{0.3\textwidth}
  \centering
  \vspace{-5pt}
  \includegraphics[width=0.28\textwidth]{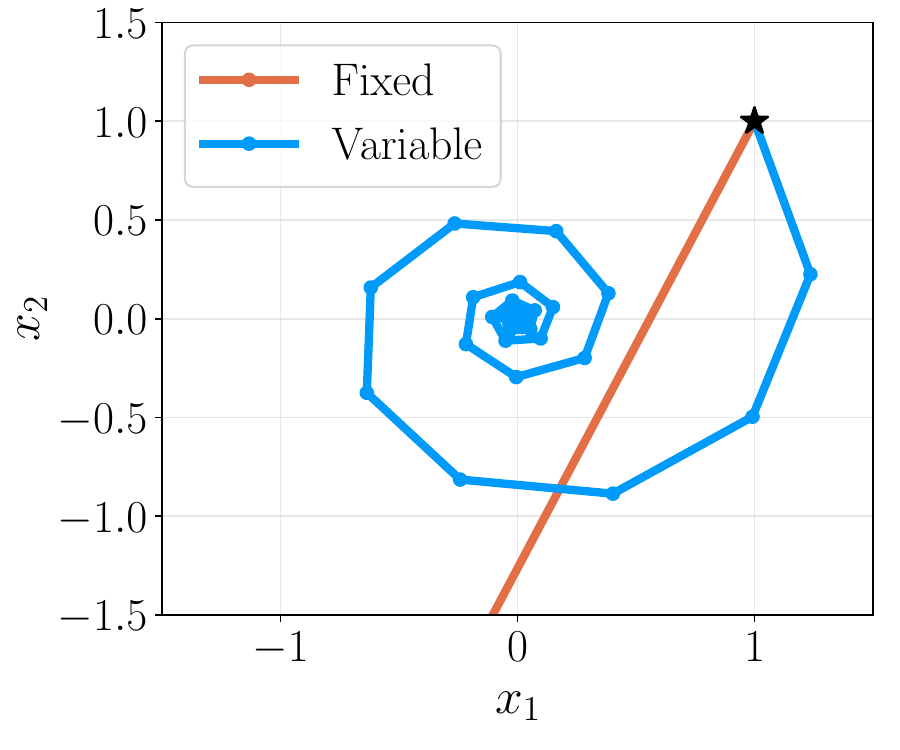}
  \caption{\eqref{eq:EG+-variable} with fixed and variable stepsizes \eqref{eq:eg+-variable-stepsizes} for the problem of \cref{ex:diverging-toilet}.}
  \label{fig:toilet}
\end{wrapfigure}

Based on \cref{prop:alpha_sym_smooth}, it is natural to select $\pre(y) = \frac{y}{L_0 + L_1 \Vert y \Vert^\alpha}$ under $\alpha$-symmetric $(L_0, L_1)$-Lipschitz continuity.
In that case, our framework yields \eqref{eq:EG+-variable} iterations with variable stepsizes
\begin{equation} \label{eq:eg+-variable-stepsizes}
    \gamma_k = \frac{\gamma}{L_0+L_1 \Vert F(x^k) \Vert^\alpha}, \qquad \omega_k = \frac{\omega}{L_0+L_1 \Vert F(z^k) \Vert^\alpha}, \qquad \gamma, \omega > 0.
\end{equation}
\Cref{fig:toilet} visualizes these adaptive stepsizes for the problem from \cref{ex:diverging-toilet} with \(\gamma = \omega = 0.95\).
The following convergence result is an immediate consequence of \cref{th:p-ceg+-convergence}. 
\begin{corollary} \label{cor:EG+-variable-monotone}
    Let $F : E \to E$ be $\alpha$-symmetric $(L_0, L_1)$-Lipschitz continuous with constants $L_0 > 0$, $L_1 \geq 0$ and $\alpha \in (0, 1]$, and suppose that \cref{assump:generalized-wmvi} holds with $C=E$, $\pre(y)=y/(L_0+L_1\|y\|^\alpha)$, and $\rho = 0$ (e.g., if $F$ is monotone). Then, the iterates of \eqref{eq:EG+-variable} with variable stepsizes given by \eqref{eq:eg+-variable-stepsizes} and $\gamma \in (0, 1), \omega \in (0, \gamma)$ converge in the sense that $\lim_{k \to \infty} \Vert F x^k \Vert = 0$.
    Moreover, for every $x^\star\in S^\star$, defining $\beta \coloneq \nicefrac{\omega}{\gamma}$ and $\bar L \coloneq L_1 \sup_{k\in\bN} \Vert F x^k \Vert^\alpha < \infty$, we have that
    \begin{equation}
        \min_{0 \leq k \leq K} \Vert F x^k \Vert^2 \leq \frac{(L_0 + \bar L)^2 \Vert x^0 - x^\star \Vert^2}{\beta \gamma^2 (1 - \beta) (1 - \gamma)^2(K+1)}.
    \end{equation}
\end{corollary}
Note that if $F$ is Lipschitz continuous with constant $L_0$, i.e., $L_1 = 0$, then $\bar L = 0$.

\begin{remark} \label{rem:big_stepsize}
    In fact, we can show convergence in \cref{cor:EG+-variable-monotone} for $0<\omega<\frac{2\gamma}{1+\gamma}$ by sharpening the proof of \cite[Corollary 3.2]{pethick2023escaping}. Specifically, since $\pre\circ F$ is $1$-Lipschitz, the map $H=\mathrm{id}-\gamma(\pre\circ F)$ is $\frac{1}{1+\gamma}$-cocoercive. Using this constant in place of the $\frac12$-cocoercivity bound in \cite[Eq.\,(B.6)]{pethick2023escaping} yields the enlarged stepsize range.
\end{remark}
Our stepsize scheme differs from those analyzed under $\alpha$-symmetric $(L_0, L_1)$-Lipschitz continuity in \cite{vankov2024generalized,choudhury2025extragradient} in three respects.
First, both \cite{vankov2024generalized,choudhury2025extragradient} propose stepsizes $\gamma_k = \omega_k$ that depend solely on $\Vert Fx^k \Vert$.
In contrast, our nonlinear preconditioning framework yields stepsizes where $\gamma_k$ depends on $\Vert F x^k \Vert$ and $\omega_k$ depends on $\Vert F z^k \Vert$.
Second, \cref{cor:EG+-variable-monotone} unifies the cases $\alpha = 1$ and $\alpha < 1$, whereas \cite{choudhury2025extragradient} requires qualitatively different stepsizes for the two cases.
Finally, \cref{cor:EG+-variable-monotone} enables the \emph{enlarged stepsize range} $0 < \gamma < 1$ compared to \cite{choudhury2025extragradient}, where the analysis considers $\gamma \leq 0.45$ for $\alpha = 1$ and significantly smaller $\gamma > 0$ if $\alpha < 1$.
For example, under $L_0$-Lipschitz continuity and with $\alpha < 1$, \cite{choudhury2025extragradient} corresponds to $\gamma \leq \frac{L_0}{2\sqrt{2}K_0} = \frac{1}{2 \sqrt{2} (2^{\alpha^2 / (1-\alpha)} + 1)} \leq \frac{1}{4 \sqrt{2}} \approx 0.177$.
A summary of supported stepsizes for \eqref{eq:EG+-variable} under $\alpha$-symmetric $(L_0, L_1)$-Lipschitz continuity is given in \cref{tab:stepsizes-comparison}.

The result from \cref{cor:EG+-variable-monotone} can be further improved from a sublinear to a linear rate by assuming $\mu$-strong monotonicity of $F$, rather than merely monotonicity. In that case, the main inequality replacing \eqref{eq:wmvi} is $\langle \pre(F x), x - x^\star \rangle \geq \mu \isocoe(\|Fx\|)\|x - x^\star\|^2$ and we obtain the following result.
\begin{theorem}\label{th:p-eg-strongly-monotone}
    Suppose that $F : E \to E$ is $\mu$-strongly monotone with $\mu>0$, $C = E$, and that \cref{assum:isotropy,def:aniso_cont} hold.
    Let $\gamma\in(0,\nicefrac{1}{L})$, $\omega\in(0,\gamma)$, $\beta = \nicefrac{\omega}{\gamma}$ and $x^\star \in \zer F$.
    Then, the iterates of \eqref{eq:EG+} converge in the sense that $\lim_{k \to \infty} \Vert F x^k \Vert = 0$. If, moreover, $\kappa$ is nonincreasing, then
    \[
        \Vert x^{k+1} - x^\star \Vert^2 \leq \left( 1 - \theta \right)^{k+1} \Vert x^0 - x^\star \Vert^2, \quad \text{where} \quad
        \theta\coloneq \tfrac{2\omega \mu \kappa(M)(1-\gamma^2L^2)}{(1-\gamma^2L^2) + 2 \gamma \mu\kappa(M)} \in (0, 1),
    \]
    with $M\coloneq \sup_{k\in\bN} \|Fz^k\|<\infty$.
\end{theorem}
This can again be specialized for $\alpha$-symmetric $(L_0, L_1)$-Lipschitz continuous operators $F$.
\begin{corollary}\label{cor:p-eg-strongly-monotone}
    Suppose that $F : E \to E$ is $\mu$-strongly monotone with $\mu>0$, $C=E$ and $F$ is $\alpha$-symmetric $(L_0, L_1)$-Lipschitz continuous with constants $L_0 > 0$, $L_1 \geq 0$ and $\alpha \in (0, 1]$.
    Then, the iterates of \eqref{eq:EG+-variable} with variable stepsizes given by \eqref{eq:eg+-variable-stepsizes} and $\gamma \in (0, 1), \omega \in (0, \gamma)$ converge in the sense that $\lim_{k \to \infty} \Vert F x^k \Vert = 0$.
    Moreover, for any $x^\star\in\zer F$, defining $\bar L \coloneq L_1\sup_{k\in\bN} \|Fz^k\|^\alpha < \infty$, we have
    \[
    \Vert x^{k+1} - x^\star \Vert^2 \leq \left( 1 - \theta \right)^{k+1} \Vert x^0 - x^\star \Vert^2, \quad \text{with} \quad
        \theta \coloneq \tfrac{2 \omega \mu (1 - \gamma^2)}{(L_0+ \bar L)  (1- \gamma^2)+ 2 \gamma \mu} \in (0, 1).
    \]
\end{corollary}

\subsection{The stochastic setting}
In this section we consider a stochastic extension of the preconditioned Extragradient method for solving \eqref{eq:rfp}. To that aim we assume access to a stochastic oracle \(F_\xi\) such that, for every $x\in E$,
\begin{equation} \label{eq:variance}
    \bE[F_\xi x] = Fx, \qquad \bE[\|F_\xi x-Fx\|^2] \leq \sigma^2.
\end{equation}
At every iteration $k\in\bN$, we draw two mutually independent minibatches of i.i.d.\ samples $\xi^k=(\xi_1^k,\ldots,\xi_{N_x}^k)$ and $\zeta^k=(\zeta_1^k,\ldots,\zeta_{N_z}^k)$, independently of all previously drawn samples, and define $\widetilde Fx^k \coloneq \tfrac{1}{N_x}\sum_{i=1}^{N_x}F_{\xi_i^k}x^k$ and $\widetilde Fz^k \coloneq \tfrac{1}{N_z}\sum_{i=1}^{N_z}F_{\zeta_i^k}z^k$.

We consider the two-timescale version of (stochastic) preconditioned Extragradient:
\begin{equation}\label{eq:EG+-variable_stoch}
    z^k = x^k - \gamma \pre(\widetilde F x^k), \qquad
    x^{k+1} = x^k - \omega \pre(\widetilde F z^k).
\end{equation}
Define the filtration corresponding to our sampling scheme by $\cF_{-1}\coloneq\sigma(\emptyset)$,
$\cG_k\coloneq\cF_{k-1}\vee\sigma(\xi^k)$ and
$\cF_k\coloneq\cG_k\vee\sigma(\zeta^k)$. By construction, $x^k$ is $\cF_{k-1}$-measurable and $z^k$ is $\cG_k$-measurable.

Although the minibatch estimators $\widetilde Fx^k$ and $\widetilde Fz^k$ are conditionally unbiased, the nonlinear preconditioning generally introduces bias, i.e., $\bE[\pre(\widetilde Fx^k) \mid \cF_{k-1}]$ need not equal $\pre(Fx^k)$. Thus, standard stochastic Extragradient analyses based on an unbiased oracle for $\pre\circ F$ cannot be applied directly. Instead, we utilize the Lipschitz continuity of $\pre$ to obtain a stochastic descent-like inequality.
\begin{theorem} \label{thm:stoch}
    Suppose that $\pre$ satisfies \cref{assum:isotropy} and is Lipschitz continuous with constant $L_\pre$. Further, suppose that $F$ satisfies \cref{def:aniso_cont} with $L > 0$ and \cref{assump:generalized-wmvi} with $\rho = 0$ and $C=E$ and let \eqref{eq:variance} hold true. Fix $x^\star\in S^\star$, and let $\{x^k\}_{k=0}^K$ be generated by \eqref{eq:EG+-variable_stoch} with $0<\gamma\leq\frac{1}{2\sqrt{1+L^2}}$, $\omega=\frac{\gamma}{\sqrt K}$, $N_x=\lceil\sqrt K\rceil$ and $N_z=K$. Then,
    \[
        \frac{1}{K}\sum_{k=0}^{K-1}\bE[\|\pre(Fx^k)\|^2] \leq \frac{4\exp(2\gamma^2)}{\gamma^2\sqrt{K}}\left [\|x^0-x^\star\|^2 + \frac{L_{\pre}^2}{2\gamma^2}\left( \frac{9}{2} \gamma^2 \sigma^2 + \sigma^2+ 4\frac{\gamma^2 \sigma^2}{K^{1/2}}\right) \right].
    \]
\end{theorem}
The large batch sizes in \cref{thm:stoch} control the bias introduced by nonlinear preconditioning. Similar large-batch requirements also appear in the VI literature, even for methods that do not use nonlinear preconditioning, see  for example \cite{kotsalis2022simple,kalogiannis2025solving}. To the best of our knowledge, \cref{thm:stoch} is the first stochastic extragradient result under the broader $\mathcal{K}$-Lipschitz and Minty condition. The closest work, \cite{vankov2024generalized2}, treats the narrower $\alpha$-symmetric class and assumes $p$-quasi-sharpness, $\langle Fx,x-x^\star \rangle \geq \mu \text{dist}^p(x,X^\star)$ for $\mu > 0$, a strong error-bound condition.

\paragraph{Experiments} In \cref{app:experiments}, we provide two synthetic experiments and two experiments for training GANs that show that preconditioning can improve performance considerably in practice. In general, preconditioning will not always improve performance when compared to nonpreconditioned algorithms (cf.\ no free lunch). Nevertheless, in many cases, a normalization-like preconditioner reduces sensitivity to initialization, since the normalization prevents excessively large updates far from the solution, allowing the method to remain stable from initial points for which a fixed-step Euclidean method would require severe tuning or may diverge.

\section{Continuous-time analysis and optimal control perspective}
\label{sec:OC}
In this section we study a continuous-time formulation of \eqref{eq:GD} as follows: find an \emph{absolutely continuous} curve $x \in X := \text{AC}([0,T],E)$ such that
\begin{equation} \label{eq:cont_time}
    \dot x(t) = -\nabla \phi^*(Fx(t)) \qquad \text{ a.e.\ on } (0,T)
\end{equation}
and $x(0) = x_0$. %
Under Lipschitz continuity of $\nabla\phi^*\circ F$ (cf.\ \cref{def:aniso_cont}), existence of a solution follows from standard arguments \cite[Corollary 2.6]{teschl2012ordinary}. The optimization properties of this flow follow in a similar manner to our analysis in the discrete-time setting. We are mainly interested in obtaining deeper insights on the proposed method by establishing a \emph{variational} perspective on the continuous-time algorithm, i.e., establishing the algorithm as a \emph{solution to an optimization problem}.

Our starting point is the Brezis--Ekeland principle \cite{brezis1976principe}, which characterizes the standard flow $\dot x(t) \in -\partial f(x(t))$ of a convex function $f$ through a minimization problem. This principle was extended to maximally monotone operators, i.e., to inclusions of the form $\dot x(t) \in -Fx(t)$, in \cite{visintin2008extension} using the \emph{Fitzpatrick} function of $F$:
\[
    f_F(x,y) \coloneq \sup_{(\bar x,\bar y) \in \gph F}\{\langle \bar y, x-\bar x \rangle + \langle y,\bar x \rangle\}.
\]
This approach covers standard flows for monotone VIs, but it does not apply directly to \eqref{eq:cont_time} because $\nabla \phi^* \circ F$ need not be monotone. The works on the so-called \emph{doubly nonlinear equations} \cite{stefanelli2008brezis} cover similar equations, but due to the fact that $F$ does not correspond to a function, the techniques used in that literature do not carry over to our setting immediately. Instead, we exploit the hypomonotonicity of $\nabla \phi^* \circ F$, which follows from its Lipschitz continuity, to obtain the following result.
\begin{theorem}[Variational formulation] \label{thm:var_form}
    Let $\phi^*$ satisfy \cref{assum:phistar}, and fix $x_0\in E$ and $T>0$.
    Suppose $F$ satisfies \cref{def:aniso_cont} with $\mathcal{K} = \nabla\phi^*$ and constant $L>0$. Define $B\coloneqq\nabla \phi^* \circ F + L\id$, $X_{x_0} \coloneqq \{ x \in \text{AC}([0,T],E) \mid x(0) = x_0\}$ and the functional
    \begin{equation}
        J : X_{x_0} \to [0, +\infty] : x \mapsto \int_0^T [f_B(x, Lx-\dot x) - L\|x\|^2]dt + \tfrac{1}{2}\|x(T)\|^2 - \tfrac{1}{2}\|x_0\|^2.
    \end{equation}
    Then $B$ is maximally monotone, $J(x) \geq 0$ for $x \in X_{x_0}$ and a curve $\bar x\in X_{x_0}$ solves \eqref{eq:cont_time} if and only if $J(\bar x) = 0$.
\end{theorem}
The problem described in \cref{thm:var_form} is reminiscent of a finite-time optimal control problem. Adapting the viewpoint of \cite{tzen2023variational,oikonomidis2026nonlinearly}, we can interpret the continuous-time algorithm under isotropy as an infinite-horizon stabilizing optimal control \cite[\S 8.5]{sontag2013mathematical}. To that aim, consider the simple dynamics
\begin{equation} \label{eq:control}
    \dot x(t) = u(t)
\end{equation}
and assume that $F$ is strictly monotone, denoting by $x^\star$ its unique zero. For $x_0 \in E$, let $U_{x_0}$ be the class of \emph{stabilizing controls} at $x_0$, i.e., all locally essentially bounded maps $u:[0,\infty)\to E$ such that $x(t)$ is defined for all $t \geq 0$ and $x(t) \to x^\star$. Consider the \emph{cost} $J_\infty(x_0,u):=\int_0^\infty q(x(s),u(s))ds$, with $q(x,u):=f_B(x,Lx-u)-L\|x\|^2+\langle u,x^\star \rangle$, where $B=\pre\circ F+L\id$, and the candidate \emph{Lyapunov function} $V(x):=\tfrac{1}{2}\|x-x^\star\|^2$. We then have the following:
\begin{theorem}[Optimal control interpretation] \label{thm:control}
    Let $\pre$ satisfy \cref{assum:isotropy} and let $F$ be strictly monotone and satisfy \cref{def:aniso_cont} with constant $L$. Fix $x_0\in E$. Then, the following hold:
    \begin{thmenum}
        \item \label{thm:control:lower_bound} For any stabilizing control $u \in U_{x_0}$ and all $t\geq 0$, we have
            $\int_0^t q(x(s),u(s))ds \geq V(x_0) - V(x(t))$
        and in particular $J_\infty (x_0,u) \geq V(x_0)$.
        \item \label{thm:control:optimal} The control $u(t) = -\pre(Fx(t))$ is optimal and stabilizing and we have
        \begin{equation*}
            J_\infty(x_0,u) = V(x_0) = \tfrac{1}{2}\|x_0-x^\star\|^2.
        \end{equation*}
        Moreover, $V(x)$ is a global control-Lyapunov function for the closed-loop system.
    \end{thmenum}
\end{theorem}
\Cref{thm:control} shows that, among all controls that stabilize \eqref{eq:control}, the control induced by our algorithm minimizes the cost $J_\infty$. Continuous-time algorithms for inclusion problems have received considerable attention \cite{lin2022continuous,boct2024extra,boct2025fast}. However, to the best of our knowledge, this is the first result linking these methods directly to variational and optimal control problems.

\section{Conclusion}
In this paper, we developed a dual-space preconditioning framework for variational inequalities and root-finding problems beyond classical cocoercivity and Lipschitz continuity. We introduced a nonlinear cocoercivity condition based on a preconditioner $\nabla\phi^*$ and established convergence of preconditioned forward methods for a broad class of reference functions. We also proposed a generalized Lipschitz condition encompassing $(L_0,L_1)$-type assumptions, enabling a unified analysis of preconditioned extragradient and optimistic methods. For isotropic preconditioners, these schemes reduce to simple variable-stepsize variants of established algorithms. We further extended the framework to stochastic settings and demonstrated its practical effectiveness numerically. Finally, we studied the associated continuous-time dynamics, deriving a variational characterization and an optimal control interpretation of the preconditioned flow. Overall, dual-space preconditioning provides a flexible and principled approach to designing operator methods beyond the classical setting. Future directions include more general composite and anisotropic schemes, stochastic methods that avoid large minibatches, and choices of $\phi$ tailored to specific problem classes.

\section*{Acknowledgments}
This work was supported by the Research Foundation Flanders (FWO) PhD grant 11A8T26N; Research Council KUL grant C14/24/103. The resources and services used in this work were provided by the VSC (Flemish Supercomputer Center), funded by the Research Foundation - Flanders (FWO) and the Flemish Government.

\printbibliography

\appendix

\newpage
\section{Preliminaries and helper results}
Let $A:E \rightrightarrows E$. We say that $A$ is \emph{$\mu$-monotone} for some $\mu \in \bR$ if
\begin{align*}
    \langle y-\bar y,x-\bar x \rangle \geq \mu \|x-\bar x\|^2 \qquad \forall (x,y),(\bar x,\bar y) \in \gph A.
\end{align*}
We also make distinctions based on the sign of the constant $\mu$. If $\mu < 0$ we say that $A$ is \emph{$|\mu|$-hypomonotone}, while if $\mu > 0$ we say that it is \emph{$\mu$-strongly monotone}. If $\mu = 0$, we say that $A$ is \emph{monotone}. We say that $A$ is $\rho$-comonotone for some $\rho \in \bR$ if
\begin{align*}
    \langle y-\bar y,x-\bar x \rangle \geq \rho \|y-\bar y\|^2 \qquad \forall (x,y),(\bar x,\bar y) \in \gph A.
\end{align*}
If $\rho < 0$ we say that $A$ is \emph{$|\rho|$-cohypomonotone}, while if $\rho > 0$ we say that it is \emph{$\rho$-cocoercive}. Clearly, cocoercive operators are single-valued.

The operator $A$ is said to be \emph{maximally (co)monotone} if its graph is not strictly contained in the graph of another (co)monotone operator. A monotone and continuous single-valued operator is maximally monotone in light of \cite[Corollary 20.28]{bauschke2017correction}.

\subsection{Useful inequalities}
We will also make use of the following standard inequalities. For $a, b\in E$ and $e > 0$, we have
\begin{align}
    2\langle a, b\rangle &\leq e\|a\|^2 + \tfrac{1}{e}\|b\|^2 \label{eq:fenchel_young}
    \\
    \|a + b\|^2 &\leq (1 + e)\|a\|^2 + (1 + \tfrac{1}{e}) \|b\|^2 \label{eq:young}
    \\
    -\|a+b\|^2 &\leq -\tfrac{1}{1+e}\|a\|^2 + \tfrac{1}{e}\|b\|^2 \label{eq:neg_young}
\end{align}

\section{Proofs of \cref{sec:problem_form}}
\subsection{Proof of \cref{thm:dual_relative}}
\begin{proof}
    Let $z, \bar z \in E$ and consider \eqref{eq:smooth_relative} with $x = z$, $\bar x = \bar z$:
    \begin{equation*}
        \phi^*(\nabla f(z)) - \phi^*(\nabla f(\bar z)) - \langle  \nabla f(z)- \nabla f(\bar z),\nabla \phi^*(\nabla f(\bar z)) \rangle \leq L [f(\bar z) - f(z) - \langle \nabla f(z),\bar z-z \rangle]
    \end{equation*}
    and \eqref{eq:smooth_relative} with $x = \bar z$, $\bar x = z$:
    \begin{equation*}
        \phi^*(\nabla f(\bar z)) - \phi^*(\nabla f(z)) - \langle \nabla f(\bar z)- \nabla f(z),\nabla \phi^*(\nabla f(z)) \rangle \leq L [f(z) - f(\bar z) - \langle \nabla f(\bar z),z-\bar z \rangle].
    \end{equation*}
    Adding the two inequalities gives
    \begin{equation*}
        \langle \nabla f(z)- \nabla f(\bar z), \nabla \phi^*(\nabla f(z))-\nabla \phi^*(\nabla f(\bar z))\rangle \leq L[\langle \nabla f(\bar z) - \nabla f(z), \bar z - z\rangle],
    \end{equation*}
    which is the claimed result.
\end{proof}

\subsection{Proof of \cref{thm:gd_general_results}} \label{subsec:gd_res}
\begin{proof}
    By definition of $D_{\phi^*}$, 
    \[
        \phi^*(Fx^{k+1}) = \phi^*(Fx^{k}) + \langle \nabla \phi^*(Fx^{k}),Fx^{k+1} - Fx^{k} \rangle + D_{\phi^*}(Fx^{k+1}, Fx^{k}).
    \]
    Substituting \eqref{eq:GD} and using the identity $D_{\phi^*}(u,v)+D_{\phi^*}(v,u)=\langle\nabla\phi^*(u)-\nabla\phi^*(v),u-v\rangle$ gives
    \begin{multline*}
          \phi^*(Fx^{k+1}) + D_{\phi^*}(Fx^k, Fx^{k+1}) = \phi^*(Fx^{k}) + \tfrac{1}{\gamma}\langle x^k-x^{k+1},Fx^{k+1} - Fx^{k} \rangle 
        \\
        + \langle \nabla \phi^*(Fx^{k+1})-\nabla \phi^*(Fx^{k}),Fx^{k+1}-Fx^k \rangle.
    \end{multline*}
    By $\nabla\phi^*$-cocoercivity of $F$, we obtain
    \[
        \phi^*(Fx^{k+1}) + D_{\phi^*}(Fx^k, Fx^{k+1}) \leq \phi^*(Fx^{k}) + (\tfrac{1}{\beta}-\tfrac{1}{\gamma})\langle Fx^{k+1}-Fx^k, x^{k+1}-x^k \rangle.
    \]
    Since $\gamma \leq \beta$ and $F$ is monotone, the sequence $\{\phi^*(Fx^k)\}_{k \in \bN}$ is nonincreasing. The summability of $D_{\phi^*}(Fx^k, Fx^{k+1})$ follows by summing up the inequality, while the fact that $\{Fx^k\}_{k \in \bN}$ is bounded is a consequence of $\phi^*(Fx^k) \leq \phi^*(Fx^0)$ which implies that $Fx^k$ stays in the $\phi^*(Fx^0)$-sublevel set of $\phi^*$. This in turn is bounded since the set $\argmin \phi^*$ is bounded in light of \cite[Proposition 11.13]{bauschke2017correction}. It remains to show that every limit point of $\{x^k\}_{k\in\bN}$ belongs to $\zer F$. We proceed with arguments similar to those in the proof of \cite[Theorem 4.19]{laude2023anisotropic}.

    Since $\{Fx^k\}_{k\in\bN}$ is bounded and $D_{\phi^*}(Fx^k, Fx^{k+1}) \to 0$, strict convexity and continuity of $\phi^*$ imply that $\|Fx^{k+1} - Fx^k\| \to 0$. Now extract a convergent subsequence $x^{k_j} \to \bar x$ and note that $Fx^{k_j} \to F \bar x$ by continuity. Moreover, we have $Fx^{k_j + m} \to F\bar x$ for every $m\in\bN$. Therefore, by the update \eqref{eq:GD}
    \[
        x^{k_j + m} = x^{k_j} - \gamma \sum_{l=0}^{m-1} \nabla\phi^*(Fx^{k_j + l}) \to \bar x - m\gamma \nabla\phi^*(F \bar x),
    \]
    which implies $F(\bar x - m\gamma \nabla\phi^*(F \bar x)) = F \bar x$ for all $m\in\bN$. Now take any $x^\star\in \zer F$. Since $\nabla\phi^*$-cocoercivity implies monotonicity of $F$, we have for every $m\in\bN$
    \[
        0 \leq \langle F\bar x, \bar x - m\gamma \nabla\phi^*(F\bar x) - x^\star\rangle = \langle F\bar x, \bar x - x^\star\rangle - m\gamma \langle F\bar x, \nabla\phi^*(F\bar x)\rangle.
    \]
    Since $\phi^*$ is even, $\nabla\phi^*(0)=0$, and $\nabla\phi^*$ is strictly monotone, we have $\langle F \bar x, \nabla\phi^*(F\bar x)\rangle > 0$ whenever $F \bar x \neq 0$. Letting $m\to\infty$ shows that $F\bar x$ is necessarily zero, from which we conclude that every limit point of $\{x^k\}_{k\in\bN}$ belongs to $\zer F$.
\end{proof}

\subsection{Proof of \cref{ex:aniso_cocoer}} \label{sec:aniso_cocoer}
\begin{proof}
It is clear that $F(x_1, x_2) = \begin{psmallmatrix}
    x_1 + x_1^3 \\ x_1 + x_2
\end{psmallmatrix}$ is not Lipschitz, and thus also not Euclidean cocoercive. Indeed, take $x = (t, 0)$ and $\bar x = (0, 0)$ which gives $\tfrac{\|Fx - F\bar x\|}{\|x - \bar x\|} = \sqrt{(1 + t^2)^2 + 1} \to \infty$ as $t\to\infty$.

Next, we show that $F$ is $\frac12$-$\nabla\phi^*$-cocoercive. Note that ${h^*}'(t) = \frac{t}{1+|t|}$, and both $t \mapsto {h^*}'(t)$ and $t\mapsto {h^*}'(t + t^3)$ are increasing and $1$-Lipschitz. Let $x,\bar x\in\bR^2$. It follows that
\[
    \langle Fx - F\bar x, \nabla\phi^*(Fx) - \nabla\phi^*(F\bar x)\rangle \leq (x_1 + x_1^3 - \bar x_1 - \bar x_1^3) (x_1 - \bar x_1) + (x_1 + x_2 - \bar x_1 - \bar x_2)^2.
\]
Moreover, since $t\mapsto t^3$ is increasing, we have $(x_1 + x_1^3 - \bar x_1 - \bar x_1^3)(x_1 - \bar x_1) \geq (x_1 - \bar x_1)^2$. Denoting $\Delta = \langle Fx - F\bar x, x - \bar x\rangle - \frac12 \langle Fx - F\bar x, \nabla\phi^*(Fx) - \nabla\phi^*(F\bar x)\rangle $, we therefore obtain
\begin{align*}
    \Delta &\geq \frac12 (x_1 - \bar x_1)^2 + (x_1 + x_2 - \bar x_1 - \bar x_2)(x_2 - \bar x_2) - \frac12 (x_1 + x_2 - \bar x_1 - \bar x_2)^2 = \frac12(x_2 - \bar x_2)^2.
\end{align*}
The nonnegativity of the square shows that $F$ is $\frac12$-$\nabla\phi^*$-cocoercive.

Finally, for any stepsize $\gamma > 0$, we can choose an initial point from which the Euclidean forward iteration $x^{k+1} = x^k - \gamma Fx^k$ diverges. Indeed, its first coordinate satisfies $x_1^{k+1}=x_1^k(1-\gamma-\gamma(x_1^k)^2)$. If $|x_1^0|\geq 2/\sqrt{\gamma}$, then $|x_1^{k+1}| \geq (3 + \gamma)|x_1^k|$ for every $k\in\bN$, so $|x_1^k|\to\infty$.
\end{proof}

\subsection{Proof of \cref{ex:aniso_not_fejer}}
\begin{proof}
    We show that $Fx = Ax$ with $A = \begin{psmallmatrix}
        1 & -2 \\ 2 & 1
    \end{psmallmatrix}$ is a $\beta$-$\nabla\phi^*$-cocoercive operator with $\beta=\tfrac{1}{5}$ and $\phi^*(y_1, y_2) = h^*(y_1) + h^*(y_2)$ with $h^*(y) = |y| - \log(|y| + 1)$. Let $x, \bar x \in \bR^2$. The $1$-Lipschitz continuity of $\nabla\phi^*$ and the Cauchy--Schwarz inequality give
    \begin{align*}
        \tfrac{1}{5} \langle Fx - F\bar x, \nabla \phi^*(Fx) - \nabla\phi^*(F\bar x)\rangle \leq \tfrac{1}{5}\|Fx - F\bar x\|^2 \leq \|x - \bar x\|^2
    \end{align*}
    where we used that the spectral norm of $A$ equals $\sqrt{5}$. Finally, writing $A = I + 2J$ with $J = \begin{psmallmatrix}
        0 &  -1 \\ 1 & 0
    \end{psmallmatrix}$ skew-symmetric yields that $\|x - \bar x\|^2 = \langle Fx - F\bar x, x - \bar x\rangle$, from which the desired result follows.
\end{proof}

\subsection{Proof of \cref{thm:iso_conv_GD}}
\begin{proof}
    Take $x^\star$ such that $Fx^\star=0$. We have:
    \begin{align*}
        \|x^{k+1}-x^\star\|^2 
        &= \|x^k-x^\star\|^2 - 2 \gamma \langle \nabla \phi^*(Fx^k),x^k-x^\star \rangle + \|x^{k+1}-x^k\|^2
        \\
        & = \|x^k-x^\star\|^2 - 2 \gamma \isocoe(\|Fx^k\|) \langle Fx^k,x^k-x^\star \rangle + \|x^{k+1}-x^k\|^2
        \\
        & \leq \|x^k-x^\star\|^2- 2\gamma \beta \isocoe(\|Fx^k\|) \langle \nabla \phi^*(Fx^k),Fx^k \rangle
        + \|x^{k+1}-x^k\|^2
        \\
        & = \|x^k-x^\star\|^2 + (1-\tfrac{2\beta}{\gamma})\|x^{k+1}-x^k\|^2
    \end{align*}
    where the first equality follows by the Pythagorean theorem, the second by the fact that $\nabla \phi^*(y) = \isocoe(\|y\|)y$, the inequality by substituting $x=x^\star$ and $\bar x = x^k$ in \eqref{eq:aniso_coco} and the final equality by \eqref{eq:GD} and $\nabla \phi^*(y) = \isocoe(\|y\|)y$. Since $\gamma < 2\beta$, the claimed result follows.

    Since $(h^*)'(\|y\|) = \|y\|\kappa(\|y\|) > 0$ for every $y \neq 0$, ${h^*}$ is strictly increasing on $[0,\infty)$. Therefore, $h^*(\|Fx^{k+1}\|) \leq h^*(\|Fx^k\|)$ implies $\|Fx^{k+1}\| \leq \|Fx^k\|$. Summing the descent inequality established above gives
    \begin{equation*}
        \sum_{k=0}^K \gamma (2\beta - \gamma)\|\nabla \phi^*(Fx^k)\|^2 \leq \|x^0-x^\star\|^2.
    \end{equation*}
    Since $\{\|Fx^k\|\}_{k \in \bN}$ is nonincreasing and ${h^*}'$ is increasing on $[0,+\infty)$, we have
    \begin{equation*}
        (K+1) \gamma(2\beta - \gamma) (\isocoe(\|Fx^{K}\|)\|Fx^{K}\|)^2 \leq \|x^0-x^\star\|^2
    \end{equation*}
    and the first inequality follows. If, moreover, $\isocoe$ is nonincreasing, then $\isocoe(\|Fx^k\|) \geq \isocoe(\|Fx^0\|)$ and the claimed rate follows.
\end{proof}

\subsection{Proof of \cref{ex:iso_coco}}
\begin{proof}
First, we establish some key identities regarding the operator $Fx = \|x\|^2(x + Jx)$. By the skew structure of $J$, we immediately have $\langle Fx, x\rangle = \|x\|^4$ for all $x\in E$. Further, since $x$ and $Jx$ are orthogonal, we have by the Pythagorean theorem that $\|Fx\|^2 = \|x\|^4\|x + Jx\|^2 = \|x\|^4(\|x\|^2 + \|Jx\|^2) = 2 \|x\|^6$, where we used that $J^*J = \id$ in the last equality. It follows that $F$ cannot be Euclidean cocoercive. Indeed, taking $\bar x = 0$, Euclidean cocoercivity would imply $2\beta \|x\|^2 \leq 1$ for every $x\in E$, which is a contradiction.

We now show that $F$ is $\frac{1}{6}$-$\nabla\phi^*$-cocoercive. Since $J^2 = -\id$, $F$ is bijective with $F^{-1}(y) = 2^{-2/3}\|y\|^{-2/3}(\id - J)y$ for $y\neq0$ and $F^{-1}(0) = 0$. By \eqref{eq:aniso_coco}, it suffices to prove that $F^{-1} - \frac{1}{6}\nabla\phi^*$ is monotone. We verify this by checking that the quadratic forms of $\nabla (F^{-1} - \frac{1}{6}\nabla\phi^*)$ are nonnegative, see \cite[Proposition 12.3]{rockafellar1998variational}. Fix $y\neq0$, set $u = y/\|y\|$ and let $v\in E$ be arbitrary. A direct calculation shows that
\begin{align*}
    \langle v, \nabla F^{-1}(y)v\rangle &= 2^{-2/3} \left(\|y\|^{-2/3}\|v\|^2 - \frac{2}{3} \|y\|^{-8/3} \langle y, v\rangle \langle v, (\id - J)y\rangle\right)
    \\ &= 2^{-2/3}\|y\|^{-2/3} \underbrace{\left(\|v\|^2 - \frac{2}{3}\langle u, v\rangle^2 + \frac{2}{3}\langle u, v\rangle \langle Ju, v\rangle\right)}_{\eqcolon A}.
\end{align*}
Since $u$ and $Ju$ are orthonormal, we may write $v = \langle u, v\rangle u + \langle Ju, v\rangle Ju + v_\perp$ where $v_\perp \perp u$ and $v_\perp \perp Ju$. Therefore, $\|v\|^2 = \langle u, v\rangle^2 + \langle Ju, v\rangle^2 + \|v_\perp\|^2$. Moreover, we have $\frac{2}{3}\langle u, v\rangle \langle Ju, v\rangle \geq -\frac{1}{6}\langle u, v\rangle^2 - \frac{2}{3}\langle Ju, v\rangle^2$ by \eqref{eq:fenchel_young}. It follows that
\[
   A \geq \|v\|^2 - \frac{5}{6}\langle u, v\rangle^2 - \frac{2}{3}\langle Ju, v\rangle^2  \geq \frac{1}{6}(\|v_\perp\|^2 + \langle u, v\rangle^2 + \langle Ju, v\rangle^2) = \frac{1}{6}\|v\|^2.
\]
Next, we have $\langle v, \nabla^2\phi^*(y)v\rangle = \frac{\|v\|^2}{1 + \|y\|} - \frac{\langle y, v\rangle^2}{\|y\|(1 + \|y\|)^2}$, so that we obtain
\[
    \left\langle v, \nabla \left(F^{-1} - \frac{1}{6}\nabla\phi^*\right)(y)v \right\rangle \geq \frac{1}{6} \left(2^{-2/3} \|y\|^{-2/3} - \frac{1}{1 + \|y\|}\right) \|v\|^2.
\]
Since $\frac{1+\|y\|}{\|y\|^{2/3}} \geq \frac{3}{2^{2/3}} > 2^{2/3}$, we have $2^{-2/3}\|y\|^{-2/3} \geq \frac{1}{1+\|y\|}$, and we conclude that the above quadratic form is nonnegative for all $y\neq 0$, hence $F^{-1} - \frac{1}{6}\nabla\phi^*$ is monotone along every line segment not containing the origin. That the same holds for all line segments follows from continuity. We conclude that $F$ is $\frac{1}{6}$-$\nabla\phi^*$-cocoercive.
\end{proof}

\subsection{\texorpdfstring{$\nabla \phi^*$}{∇ϕ∗}-cocoercivity of port-Hamiltonian equilibrium operators}
\label{app:port_hamiltonian}
Consider the port-Hamiltonian \cite{van2014port} equilibrium operators $Fx = (\id + \tau J)\nabla H(x)$ with $J^* J = \id$, $J^* = -J$ and the Hamiltonian $H(x) = \frac{1}{p + 1}\|x - x^\star\|^{p+1}$ with $p > 1$. These operators are used extensively in mechanical, electrical, hydraulic, and other multiphysics systems. Let $\nabla\phi^*(y) = (1 + \tau^2)^{(p-1)/(2p)}\|y\|^{1/p-1}y$ for $y\neq 0$ and $\nabla\phi^*(0) = 0$. Suppose that $|\tau| < \frac{2\sqrt{p}}{p - 1}$. Then, $F$ is $\beta$-$\nabla\phi^*$-cocoercive with $\beta = \frac{1 - \frac{|\tau|(p-1)}{2\sqrt{p}}}{1 + \tau^2} > 0$.

Indeed, let $x\neq x^\star$ and $v\in E$. Since $\nabla^2H(x) = \|x - x^\star\|^{p-1}\id + (p - 1)\|x - x^\star\|^{p-3}(x - x^\star)(x - x^\star)^*$, we have
\begin{align*}
    |\langle v, J\nabla^2 H(x)v\rangle| &= (p-1)\|x - x^\star\|^{p - 3} |\langle x - x^\star, v\rangle\langle J(x - x^\star), v\rangle| \\
    &\leq \frac{p-1}{2\sqrt{p}} \|x - x^\star\|^{p-3} \left((p - 1)\langle x - x^\star, v\rangle^2 + \|x - x^\star\|^2 \|v\|^2\right) \\
    &= \frac{p-1}{2\sqrt p} \langle v, \nabla^2H(x)v\rangle,
\end{align*}
where the first equality follows from $J^* = -J$, and the inequality follows from \eqref{eq:fenchel_young} with $e=\sqrt{p}$ and the fact that $x -x^\star$ and $J(x - x^\star)$ are orthogonal and have equal norm. The same bound holds at $x = x^\star$, since $\nabla^2H(x^\star)=0$. By the fundamental theorem of calculus, for $x,\bar x\in E$, we have
\begin{align*}
    |\langle J(\nabla H(x) - \nabla H(\bar x)), x - \bar x\rangle| &\leq \frac{p-1}{2\sqrt{p}} \int_0^1 \langle \nabla^2 H(\bar x + t(x - \bar x))(x - \bar x), x - \bar x\rangle dt \\
    &= \frac{p-1}{2\sqrt{p}} \langle \nabla H(x) - \nabla H(\bar x), x - \bar x\rangle,
\end{align*}
which implies $\langle Fx - F\bar x, x - \bar x\rangle \geq (1 - \frac{|\tau|(p - 1)}{2\sqrt{p}}) \langle \nabla H(x) - \nabla H(\bar x), x - \bar x\rangle$. Moreover, $\nabla\phi^*(Fx) - \nabla\phi^*(F\bar x) = (\id + \tau J)(x - \bar x)$, hence $\langle Fx - F\bar x, \nabla\phi^*(Fx) - \nabla\phi^*(F\bar x)\rangle = (1 + \tau^2)\langle \nabla H(x) - \nabla H(\bar x), x - \bar x\rangle$, where we used that $(\id + \tau J)^*(\id + \tau J) = (1 + \tau^2)\id$. It follows that
\[
    \langle Fx - F\bar x, x - \bar x\rangle \geq \frac{1 - \frac{|\tau|(p-1)}{2\sqrt{p}}}{1 + \tau^2} \langle Fx - F\bar x, \nabla\phi^*(Fx) - \nabla\phi^*(F\bar x)\rangle,
\]
hence $F$ is $\beta$-$\nabla\phi^*$-cocoercive.

Although this is a sharp theoretical example of $\nabla\phi^*$-cocoercivity where our proposed methods have good convergence guarantees, finding the equilibrium is trivial. A more practical case would be to consider the operator $Fx = (\id + \tau J)\nabla H(x) + b$ for some load $b$ and the constrained system $0 \in Fx + N_C(x)$ where $C$ represents some physical constraint.

\subsection{Discussion of \cref{subsec:extension}} \label{app:discuss_composite}
In the case $F\equiv 0$, the iteration \eqref{eq:for_back} reduces to $x^{k+1} = (\id + \gamma \nabla\phi^* \circ G)^{-1}(x^k)$, which is the iteration studied in \cite{laude2023anisotropic}. If $G$ is maximally monotone, then the residual $R_\gamma(x) = \nabla \phi\left(\frac{x - (\id + \gamma\nabla\phi^* \circ G)^{-1}(x)}{\gamma}\right)$ is $\gamma$-$\nabla\phi^*$-cocoercive. Indeed, take $x,\bar x\in E$ and define $x^+ = x - \gamma \nabla\phi^*(R_\gamma(x))$ and $\bar x^+ = \bar x - \gamma \nabla\phi^*(R_\gamma(\bar x))$. Then,
\begin{align*}
    \langle R_\gamma(x) - R_\gamma(\bar x), x - \bar x\rangle &= \begin{multlined}[t] \langle R_\gamma(x) - R_\gamma(\bar x), x^+ - \bar x^+\rangle \\ + \gamma \langle R_\gamma(x) - R_\gamma(\bar x), \nabla\phi^*(R_\gamma(x)) - \nabla\phi^*(R_\gamma(\bar x))\rangle \end{multlined} \\
    &\geq \gamma \langle R_\gamma(x) - R_\gamma(\bar x), \nabla\phi^*(R_\gamma(x)) - \nabla\phi^*(R_\gamma(\bar x))\rangle 
\end{align*}
where we used $R_\gamma(x) \in Gx^+$, $R_\gamma(\bar x) \in G\bar x^+$, and the monotonicity of $G$ in the inequality.

\section{Proofs of \Cref{sec:prec_lip}}
\subsection{Proof of \Cref{prop:alpha_sym_smooth}} \label{subsec:proof_l0l1}
We first show a preliminary result that will be used in the proof.
\begin{lemma} \label{lem:prec_bound}
    Let $y,\bar y \in E$, $L_0 > 0$, $L_1\geq0$, and $\alpha\in(0,1]$. Then, 
    \[
        \left\|\frac{y}{L_0 + L_1\|y\|^\alpha} - \frac{\bar y}{L_0 + L_1\|\bar y\|^\alpha}\right\| \leq \frac{\|y - \bar y\|}{\sqrt{(L_0 + L_1\|y\|^\alpha)(L_0 + L_1\|\bar y\|^\alpha) }}.
    \]  
\end{lemma}
\begin{proof}
    Let $y, \bar y \in E$. Without loss of generality, assume $\|y\|\geq\|\bar y\|$. Denote $Y = L_0 + L_1\|y\|^\alpha$ and $\bar Y = L_0 + L_1\|\bar y\|^\alpha$. By simple algebra, we have
     \begin{align*}
         \|\bar Y y - Y \bar y\|^2 &= \bar Y^2 \|y\|^2 + Y^2\|\bar y\|^2 - 2 Y\bar Y\langle y, \bar y\rangle \\
         &= \bar Y^2 \|y\|^2 + Y^2\|\bar y\|^2 + Y\bar Y(\|y - \bar y\|^2 - \|y\|^2 - \|\bar y\|^2) \\
         &= Y\bar Y\|y - \bar y\|^2 - (Y - \bar Y) (\bar Y\|y\|^2 - Y \|\bar y\|^2).
     \end{align*}
    Observe that
    \begin{align*}
        \bar Y\|y\|^2 - Y\|\bar y\|^2 &= (L_0 + L_1\|\bar y\|^\alpha)\|y\|^2 - (L_0 + L_1\|y\|^\alpha)\|\bar y\|^2 \\
        &= L_0(\|y\|^2 - \|\bar y\|^2) + L_1\|y\|^2\|\bar y\|^\alpha - L_1\|y\|^\alpha\|\bar y\|^2 \\
        &= L_0(\|y\|^2 - \|\bar y\|^2) + L_1\|y\|^\alpha\|\bar y\|^\alpha(\|y\|^{2-\alpha} - \|\bar y\|^{2-\alpha})
    \end{align*}
    is nonnegative since $\|y\|\geq \|\bar y\|$ and $\alpha\in(0,1]$. Therefore, $(Y-\bar Y)(\bar Y\|y\|^2 - Y\|\bar y\|^2) \geq 0$ and hence $\|\bar Yy - Y\bar y\|^2 \leq Y\bar Y\|y - \bar y\|^2$. Taking square roots and dividing through by $Y\bar Y$ (which is nonzero by our hypotheses), we conclude $ \left\|\frac{y}{Y} - \frac{\bar y}{\bar Y}\right\| \leq \frac{\|y - \bar y\|}{\sqrt{Y\bar Y}}$ which is the claimed result.
\end{proof}

We now prove \cref{prop:alpha_sym_smooth}.
\begin{proof}
    Let $x,\bar x\in E$. If $x=\bar x$, then \eqref{eq:aniso_cont} is clear. Assume that $x\neq \bar x$ and parameterize the line segment between these points as $\gamma(t) = (1-t)x + t\bar x$ for $t\in[0,1]$. By continuity of $F$, we have that $D = [0, \max_{t\in[0,1]} \|F\gamma(t)\|]$ is compact. Define
    \[
        R(\beta,\omega) = \frac{L_0 + L_1(\beta + \omega)^\alpha}{\sqrt{(L_0 + L_1\max(0, \beta - \omega)^\alpha)(L_0 + L_1\beta^\alpha) }}
    \]
    which converges uniformly to $1$ for all $\beta\in D$ as $\omega\downarrow0$ by compactness of $D$. In particular, let $\varepsilon > 0$ be arbitrary and choose $\omega > 0$ such that $R(\beta,\omega) \leq 1+\varepsilon$ for all $\beta\in D$. 
    
    For this fixed $\omega$, we invoke uniform continuity of $F\circ \gamma$ to obtain an integer $n>0$ such that on each subinterval $[i/n, (i+1)/n]$ with $i=0, \ldots, n-1$, the variation of $\|F\circ\gamma\|$ is bounded by $\omega$. Let $x_i = \gamma(i/n)$ denote the endpoints of these intervals so that we obtain
    \[
        \max_{t\in[i/n, (i+1)/n]} \|F\gamma(t)\| \leq \|Fx_i\| + \omega \quad \textnormal{and} \quad \|Fx_{i+1}\| \geq \max(0, \|Fx_i\| - \omega).
    \]
    Applying \eqref{eq:alpha_sym_smooth} on some segment $[x_i, x_{i+1}]$ then yields
    \[
        \|Fx_{i+1} - Fx_i\| \leq (L_0 + L_1(\|Fx_i\|+\omega)^\alpha) \|x_{i+1} - x_i\|.
    \]
    By \Cref{lem:prec_bound}, we then obtain
    \begin{align*}
        \|\pre(Fx_{i+1}) - \pre(Fx_i)\| &\leq \frac{\|Fx_{i+1} - Fx_i\|}{\sqrt{(L_0 + L_1\|Fx_{i+1}\|^\alpha)(L_0 + L_1\|Fx_i\|^\alpha)}} \\
        &\leq \frac{(L_0 + L_1(\|Fx_i\| + \omega)^\alpha)\|x_{i+1} - x_i\|}{\sqrt{(L_0+L_1\max(0, \|Fx_i\| - \omega)^\alpha)(L_0 + L_1\|Fx_i\|^\alpha)}} \\
        &= R(\|Fx_i\|, \omega) \frac{1}{n}\|x - \bar x\|.
    \end{align*}
    Since $\|Fx_i\| \in D$, we have $R(\|Fx_i\|, \omega)\leq 1 + \varepsilon$. Summing this bound over all $n$ subintervals and using the triangle inequality, we conclude
    \[
        \|\pre(Fx) - \pre(F\bar x)\| \leq \sum_{i=0}^{n-1} \|\pre(Fx_{i+1}) - \pre(Fx_i)\| \leq \sum_{i=0}^{n-1} (1+\varepsilon) \frac{1}{n} \|x - \bar x\| = (1+ \varepsilon) \|x - \bar x\|.
    \]
    As $\varepsilon > 0$ is arbitrary, the desired result follows.
\end{proof}

\subsection{Proof of \Cref{ex:diverging-toilet}} \label{sec:prf-diverging-toilet}
\begin{proof}
We first prove the $\pre$-Lipschitz continuity claim. Let $J = \begin{psmallmatrix}
    0 & -1 \\ 1 & 0
\end{psmallmatrix}$, $x,\bar x\in \bR^2$ and denote $G : x\mapsto \frac{e^{\|x\|^2/2}}{1+\|x\|e^{\|x\|^2/2}}x$. Then,
\[\|\pre(Fx) - \pre(F\bar x)\| = \|J(Gx - G\bar x)\| = \|Gx - G\bar x\|\]
where we used that $\|Jz\| = \|z\|$ for all $z\in\bR^2$. It therefore suffices to show that $G$ is $1$-Lipschitz. For $x\neq 0$, and denoting $r = \|x\|$, a calculation shows that
\[
    \nabla G(x)= \frac{e^{r^2/2}}{1 + re^{r^2/2}} I + \frac{e^{r^2/2} (r - e^{r^2/2})}{r(1+ re^{r^2/2})^2} xx^\top = \frac{e^{r^2/2}}{1 + re^{r^2/2}} P_x^\perp + \frac{e^{r^2/2}(1+r^2)}{(1+ re^{r^2/2})^2} P_x
\]
where $P_x$ denotes the projector $\frac{xx^\top}{\|x\|^2}$ and $P_x^\perp = I - P_x$. It follows that the operator norm equals
\[
    \|\nabla G(x)\|_{\rm op} = \max\left\{\frac{e^{r^2/2}}{1+r e^{r^2/2}}, \frac{e^{r^2/2}(1+r^2)}{\bigl(1+r e^{r^2/2}\bigr)^2} \right\}.
\]
Both quantities are bounded above by $1$ for $r>0$, hence $\|\nabla G(x)\|_{\rm op} \leq 1$ for $x\neq 0$. For $x=0$, we have $\nabla G(0)=I$, so the same bound holds. We conclude that $G$ is $1$-Lipschitz, hence $F$ is $\pre$-Lipschitz continuous with constant $1$.

It remains to show that $F$ is not $\alpha$-symmetric $(L_0,L_1)$-Lipschitz. Suppose otherwise. Then, for any $u \in \bR^2$ with $\|u\|=1$, we have $\|\nabla F(x) u\| \leq L_0 + L_1 \|Fx\|^\alpha$. Taking $x = te_1$ with $t>0$ and $u = e_1$, we have
\[
    \|\nabla F(te_1)e_1\| = \|e^{t^2/2}(Je_1 + \langle te_1, e_1\rangle J(te_1))\|  = e^{t^2/2}(1+ t^2), \qquad \|F(te_1)\| = te^{t^2/2}
\]
which implies that
\[
   e^{t^2/2}(1+t^2) \leq L_0 + L_1 t^\alpha e^{\alpha t^2/2} \iff 1 \leq \frac{L_0 e^{-t^2/2}}{1+t^2} + \frac{L_1t^\alpha e^{-(1-\alpha)t^2/2}}{1 + t^2}.
\]
The right-hand side of the latter inequality goes to $0$ as $t\to\infty$ for every $\alpha\in(0, 1]$, which is a contradiction.

\end{proof}

\subsection{Proof of \Cref{rem:alpha_sym_conv}} \label{sec:alpha_sym_conv}
\begin{proof}
    Let $y, \bar y\in E$, $L_0 > 0, L_1 \geq 0$ and $\alpha\in(0,1)$. We show that
    \[
        (1-\alpha)\|y - \bar y\| \leq (L_0 + L_1 \max(\|y\|, \|\bar y\|)^{\alpha}) \left\|\frac{y}{L_0 + L_1 \|y\|^\alpha} - \frac{\bar y}{L_0 + L_1\|\bar y\|^\alpha}\right\|
    \]
    from which the desired result immediately follows. If either $y=0$ or $\bar y=0$, then this inequality is clear.

    Without loss of generality, assume that $\|y\| \geq \|\bar y\| > 0$. Let $Y = L_0 + L_1\|y\|^\alpha$ and $\bar Y = L_0 + L_1 \|\bar y\|^\alpha$. Since $y - \bar y = Y(\tfrac{y}{Y} - \tfrac{\bar y}{\bar Y}) + \frac{Y - \bar Y}{\bar Y}\bar y$, we have from the triangle inequality that
    \[
        \|y -\bar y\| \leq Y \left\|\frac{y}{Y} - \frac{\bar y}{\bar Y}\right\| + \frac{L_1(\|y\|^\alpha - \|\bar y\|^\alpha)}{\bar Y} \|\bar y\| \leq Y \left\|\frac{y}{Y} - \frac{\bar y}{\bar Y}\right\| + (\|y\|^\alpha - \|\bar y\|^\alpha) \|\bar y\|^{1-\alpha}
    \]
    where we used that $\tfrac{L_1\|\bar y\|^\alpha}{L_0 + L_1\|\bar y\|^\alpha} \leq 1$ in the second inequality. Next, by concavity of the function $f(x) = x^\alpha$ for $x \geq 0$, we have $\|y\|^\alpha \leq \|\bar y\|^\alpha + \alpha \|\bar y\|^{\alpha -1} (\|y\| - \|\bar y\|) \iff (\|y\|^\alpha - \|\bar y\|^\alpha) \|\bar y\|^{1-\alpha} \leq \alpha (\|y\| - \|\bar y\|)$. Then, using the reverse triangle inequality, we obtain the desired result
    \[
         \|y -\bar y\| \leq Y \left\|\frac{y}{Y} - \frac{\bar y}{\bar Y}\right\| + \alpha \|y - \bar y\| \iff (1-\alpha)\|y-\bar y\| \leq Y\left\| \frac{y}{Y} - \frac{\bar y}{\bar Y}\right\|.
    \]
\end{proof}

\subsection{Proof of \Cref{lem:zero-correspondance}}
\begin{proof}
    For any $x\in C$, we have
    \[
        0 \in \pre(Fx) + N_C(x) \iff -\kappa(\|Fx\|) Fx \in N_C(x) \iff -Fx \in N_C(x) \iff 0 \in Fx + N_C(x)
    \]
    where the second equivalence follows because $N_C(x)$ is a cone, see \cite[Proposition 6.5]{rockafellar1998variational}, and $\kappa(\|Fx\|) > 0$. We conclude that $\zer(\pre \circ F + N_C) = \zer(F + N_C)$.
\end{proof}

\subsection{Proof of \Cref{lem:T-wmvi}}
\begin{proof}
    Let $x^\star \in \zer(F + N_C)$ and $(x, v) \in \gph (\pre \circ F + N_C)$. Then, $v = \kappa(\|Fx\|) Fx + u$ for some $u\in N_C(x)$ and $-Fx^\star \in N_C(x^\star)$. By the definition of the normal cone, we have $\langle Fx^\star, x - x^\star\rangle \geq 0$ and $\langle u, x - x^\star\rangle\geq 0$. Moreover, monotonicity of $F$ gives $\langle Fx - Fx^\star, x-x^\star\rangle\geq0$, hence $\langle Fx, x - x^\star\rangle \geq \langle Fx^\star, x - x^\star\rangle \geq 0$. Therefore, utilizing $\kappa(\|Fx\|) > 0$, we have
    \[
        \langle v, x - x^\star\rangle = \kappa(\|Fx\|) \langle Fx, x - x^\star\rangle + \langle u, x - x^\star \rangle \geq 0
    \]
    so the weak MVI holds with $\rho = 0$ and $S^\star = \zer(F + N_C) = \zer(\pre \circ F + N_C)$, where the latter equality is due to \cref{lem:zero-correspondance}.
\end{proof}

\subsection{Proof of \Cref{cor:EG+-variable-monotone}}
\begin{proof}
    For $C=E$, the updates \eqref{eq:EG+} reduce to
    \begin{equation*}
        z^k = x^k - \gamma \pre(F x^k), \qquad x^{k+1} = x^k - \omega \pre(F z^k),
    \end{equation*}
    where $\omega=\beta\gamma$. Take $\pre : y \mapsto \frac{y}{L_0 + L_1\|y\|^\alpha}$, so that $\pre \circ F$ is $1$-Lipschitz by \cref{prop:alpha_sym_smooth}. Moreover, $\rho = 0$ by the hypothesis, hence \cref{th:p-ceg+-convergence} applies for $0 < \gamma < 1$ and $0 < \beta < 1$, or equivalently, $0 < \omega < \gamma < 1$.
    
    From \(
    z^k - x^k = - \gamma \pre(Fx^k)
    \)
    we find
    \begin{equation*}
        \Vert Hz^k - H x^k \Vert^2 \geq (1 - \gamma L)^2 \| z^k - x^k \|^2 = (1 - \gamma L)^2 \gamma^2 \| \pre(F x^k) \|^2,
    \end{equation*}
    where the inequality follows from \cite[Lemma A.3 (ii)]{pethick2023escaping}. Applying this lower bound to the result from \cref{th:p-ceg+-convergence} with $L = 1$ yields
    \begin{equation*}
        \min_{0 \leq k \leq K} \Vert \pre(F x^k ) \Vert^2 \leq \frac{\Vert x^0 - x^\star \Vert^2}{\beta(1-\beta) \gamma^2 (1 - \gamma)^2(K+1)}
    \end{equation*}
    and that $\{z^k\}_{k\in\bN}$ is bounded. Invoking \cite[Lemma A.3 (ii)]{pethick2023escaping} again, we obtain $\|z^k- x^\star\| = \|Hx^k - Hx^\star\| \geq (1-\gamma)\|x^k - x^\star\|$, which implies that $\{x^k\}_{k\in\bN}$ is bounded. Continuity of $F$ then ensures that $\bar L = L_1\sup_{k\in\bN}\|Fx^k\|^\alpha < \infty$. Moreover, by \cref{th:p-ceg+-convergence}, every limit point of $\{z^k\}_{k\in\bN}$ belongs to $\zer F$. It follows that $\pre(Fz^k) \to 0$. Furthermore, since $\pre\circ F$ is $1$-Lipschitz, we have $\|\pre (Fz^k)\| \geq \|\pre(Fx^k)\| - \|\pre(Fz^k) - \pre(Fx^k)\| \geq (1-\gamma)\|\pre(Fx^k)\|$ where we used the update rule and $\gamma<1$. The desired result now follows from
    \[
        \Vert \pre(F x^k) \Vert = \frac{\Vert F x^k \Vert}{L_0 + L_1 \Vert F x^k \Vert^\alpha} \geq \frac{\Vert F x^k  \Vert}{L_0 + \bar L}.
    \]
\end{proof}

\subsection{Proof of \Cref{th:p-eg-strongly-monotone}}
We first provide a preliminary lemma in a similar spirit to \cref{lem:T-wmvi}.
\begin{lemma}\label{lem:eg-strong-monotonicity}
    Let $\pre$ satisfy \cref{assum:isotropy}.
    Let $F$ be $\mu$-strongly monotone with $\mu > 0$.
    Then, $\zer(\pre\circ F) = \zer F$.
    Moreover, for all $x^\star \in \zer (\pre\circ F)$ we have
    \begin{equation}
        \langle z, x - x^\star \rangle \geq \mu \isocoe(\Vert F x \Vert) \Vert x - x^\star \Vert^2, \qquad \forall (x, z) \in \gph (\pre\circ F).
    \end{equation}
\end{lemma}
\begin{proof}
   That $\zer (\pre\circ F) = \zer F$ follows from $\pre(y) = 0 \iff y = 0$. As for the second part, let $x^\star \in \zer (\pre \circ F) = \zer F$. Then, for any $(x, z) \in \gph (\pre \circ F)$ we have
    \[
        \langle z, x - x^\star \rangle = \isocoe(\Vert Fx \Vert) \langle Fx, x - x^\star \rangle = \kappa(\|Fx\|) \langle Fx - Fx^\star, x - x^\star\rangle \geq \mu \isocoe(\Vert F x \Vert) \Vert x - x^\star \Vert^2,
    \]
    where the inequality follows from $\mu$-strong monotonicity of $F$.
\end{proof}

We now move on to the proof of \Cref{th:p-eg-strongly-monotone}.
\begin{proof}
    In this setting, the updates of \eqref{eq:EG+} reduce to
    \begin{equation}
        z^k = x^k - \gamma \pre(F x^k), \qquad x^{k+1} = x^k - \omega \pre(F z^k).
    \end{equation}
    For ease of notation, we denote $T \coloneq \pre \circ F$ in the remainder of this proof.
    By the update rule and \Cref{lem:eg-strong-monotonicity}, we have
    \begin{align*}
        \|x^{k+1} - x^\star\|^2 &= \|x^k - x^\star\|^2 -2 \omega \langle Tz^k, x^k - x^\star\rangle + \omega^2\|Tz^k\|^2 \\
        &= \|x^k - x^\star\|^2 -2 \omega \langle Tz^k, z^k - x^\star\rangle - 2\gamma \omega \langle T z^k, Tx^k\rangle + \omega^2\|Tz^k\|^2 \\
        &\leq \|x^k - x^\star\|^2 -2\omega\mu\isocoe(\|Fz^k\|)\|z^k - x^\star\|^2 - 2\gamma \omega \langle T z^k, Tx^k\rangle + \omega^2\|Tz^k\|^2.
    \end{align*}
    
    Now, using the identity $-2\langle a, b\rangle = \|a - b\|^2 - \|a\|^2 - \|b\|^2$ with $a = Tz^k$ and $b=Tx^k$ yields
    \begin{align*}
        \|x^{k+1} - x^\star\|^2 &\leq \|x^k - x^\star\|^2 -2\omega\mu\isocoe(\|Fz^k\|)\|z^k - x^\star\|^2 + \gamma \omega \|Tz^k - Tx^k\|^2  \\ & \mspace{280mu} - \gamma \omega \|Tz^k\|^2 - \gamma \omega\|Tx^k\|^2
         + \omega^2\|Tz^k\|^2 \\
        &\leq  \|x^k - x^\star\|^2 -2\omega\mu\isocoe(\|Fz^k\|)\|z^k - x^\star\|^2 + \gamma \omega \|Tz^k - Tx^k\|^2  - \gamma \omega \|Tx^k\|^2
    \end{align*}
    where we dropped the nonpositive term $(\omega^2 - \gamma \omega) \|Tz^k\|^2$ in the last inequality.

    Then, by preconditioned Lipschitz continuity, \Cref{def:aniso_cont}, and the update rule, we obtain
    \begin{align*}
        \|x^{k+1} - x^\star\|^2 &\leq  \|x^k - x^\star\|^2 -2\omega\mu\isocoe(\|Fz^k\|)\|z^k - x^\star\|^2 + \gamma \omega L^2 \|z^k - x^k\|^2  - \gamma \omega \|Tx^k\|^2 \\
        &\leq  \|x^k - x^\star\|^2 -2\omega\mu\isocoe(\|Fz^k\|)\|z^k - x^\star\|^2 - \left(\frac{\omega}{\gamma} - \gamma \omega L^2 \right)\|\gamma Tx^k\|^2 \\
        &= \|x^k - x^\star\|^2 - A_k\|z^k - x^\star\|^2 - B\|\gamma Tx^k\|^2,
    \end{align*}
    where we define $A_k \coloneqq 2\omega\mu\isocoe(\|Fz^k\|)$ and $B \coloneqq \frac{\omega}{\gamma} - \gamma \omega L^2$. Since $\isocoe(\|Fz^k\|) > 0$, we have $A_k > 0$.
    Note also that by our parameter choices we have $0 < B < 1$. Therefore $\|x^{k+1} - x^\star\|^2 \leq \|x^k - x^\star\|^2$ and $\{x^k\}_{k\in\bN}$ is bounded. Summing this inequality shows that $Tx^k \to 0$, which implies that $\|Fx^k\| \to 0$.

    Suppose now that $\kappa$ is nonincreasing. Since $\{x^k\}_{k\in\bN}$ is bounded and $Tx^k \to 0$, we have that $z^k = x^k - \gamma Tx^k$ is bounded. Therefore, $M = \sup_{k\in\bN} \|Fz^k\| < \infty$. It follows that $\kappa(\|Fz^k\|) \geq \kappa(M)$, and hence $A_k \geq A\coloneqq 2\omega \mu \kappa(M) > 0$, and the descent inequality becomes
    \[
        \|x^{k+1} - x^\star\|^2 \leq \|x^k - x^\star\|^2 - A\|z^k - x^\star\|^2 - B\|\gamma Tx^k\|^2.
    \]
    We now aim to upper bound the last two terms to establish contraction. To this end, observe that Young's inequality \eqref{eq:young} with $a=z^k - x^\star$, $b=\gamma Tx^k$, and $e=A/B$ gives
    \[
        \|z^k - x^\star + \gamma Tx^k\|^2 \leq  \frac{A+B}{B} \|z^k - x^\star\|^2 + \frac{A+B}{A}\|\gamma Tx^k\|^2.
    \]
    Noting that the left-hand side is exactly $\|x^k - x^\star\|^2$ by the update rule, and multiplying both sides by $-\frac{A B}{A+B}$ yields
    \[
        - A\|z^k - x^\star\|^2 - B\|\gamma Tx^k\|^2 \leq -\frac{A B}{A+B}\|x^k - x^\star\|^2.
    \]
    Defining $\theta \coloneqq \frac{A B}{A+B}$, which satisfies $\theta \in (0, 1)$, we conclude
    \[
        \|x^{k+1} - x^\star\|^2 \leq (1 - \theta) \|x^k - x^\star\|^2,
    \]
    from which the desired result follows.
\end{proof}

\subsection{Proof of \cref{cor:p-eg-strongly-monotone}}
\begin{proof}
    This result follows from \cref{th:p-eg-strongly-monotone} by selecting the preconditioner $\pre(y) = \frac{y}{L_0 + L_1 \|y\|^\alpha}$. By \cref{prop:alpha_sym_smooth}, we have $L=1$, and $\kappa(r) = \frac{1}{L_0 + L_1 r^\alpha}$ is nonincreasing. Hence, \cref{th:p-eg-strongly-monotone} yields $\|Fx^k\| \to 0$ and the claimed linear convergence rate.

    Indeed, letting $M \coloneqq \sup_{k\in\bN}\|Fz^k\| < \infty$, we have $\bar L = L_1 M^\alpha$ and therefore $\kappa(M) = \frac{1}{L_0+\bar L}$. Substituting $L=1$ and this expression for $\kappa(M)$ into the rate constant of \cref{th:p-eg-strongly-monotone} gives
    \[
        \theta = \frac{2\omega\mu(1-\gamma^2)}
        {(L_0+\bar L)(1-\gamma^2)+2\gamma\mu},
    \]
    which proves the result.
\end{proof}

\subsection{Proof of \cref{thm:stoch}} \label{subsec:stoch}
We will require the following standard result:
\begin{lemma} \label{thm:sequences}
    Let $\{V_k\}, \{u_k\}$ be nonnegative sequences that satisfy the recursive inequality
    \begin{equation} \label{eq:rec}
        V_{k+1} \leq (1+\alpha)V_k + C - u_k
    \end{equation}
    for some $C \geq 0$ and $\alpha > 0$. If $\alpha = \tfrac{\bar\alpha}{K}$ for some $\bar\alpha >0$, then we have the following for $K\geq 1$:
    \begin{equation*}
        \sum_{k=0}^{K-1}u_k \leq \exp(\bar \alpha) (V_0 + \tfrac{C K}{\bar \alpha}).
    \end{equation*}
\end{lemma}
\begin{proof}
    Consider the rescaled quantity $W_k = (1+\alpha)^{-k}V_k$. From \eqref{eq:rec} we have:
    \begin{align*}
        W_{k+1} &\leq (1+\alpha)^{-(k+1)}[(1+\alpha) V_k + C - u_k]
        = W_k + (1+\alpha)^{-(k+1)}C - (1+\alpha)^{-(k+1)}u_k.
    \end{align*}
    Summing from $k=0$ to $k=K-1$ gives 
    \[
        \sum_{k=0}^{K-1}(1+\alpha)^{-(k+1)}u_k \leq W_0 - W_{K} + \sum_{k=0}^{K-1}(1+\alpha)^{-(k+1)}C \leq V_0 + \frac{1}{\alpha}C.
    \]
    Since $1+\alpha > 1$ we have $(1+\alpha)^{-(k+1)} \geq (1+\alpha)^{-K}$ for all $k=0,1,\ldots,K-1$, which leads us to
    \[
        \sum_{k=0}^{K-1}u_k \leq (1+\alpha)^K (V_0 + \tfrac{C}{\alpha}).
    \]
    Finally, using $\alpha=\frac{\bar\alpha}{K}$ and $\left(1+\frac{\bar\alpha}{K}\right)^K\leq \exp(\bar\alpha)$ gives the desired result.
\end{proof}
We now prove \cref{thm:stoch}. Recall that $\pre$ is an $L_{\pre}$-Lipschitz operator and the algorithm is given by
\begin{align}
    z^k &= x^k - \gamma \pre(\dOne) \label{eq:z_update}
    \\
    x^{k+1} &= x^k - \omega \pre(\dTwo) \label{eq:x_update}
\end{align}
\begin{proof}
    Fix $x^\star \in S^\star$ and use the Pythagorean theorem on $\|x^{k+1}-x^\star\|^2$:
    \begin{align*}
        \|x^{k+1}-x^\star\|^2 
        &= \|x^k -x^\star\|^2 - 2\omega \langle \pre(\dTwo),x^k-x^\star \rangle + \|x^{k+1}-x^k\|^2
        \\
        &= \|x^k -x^\star\|^2 \underbrace{- 2\omega \langle \pre(\dTwo),z^k-x^\star \rangle}_A \underbrace{- 2\omega \gamma \langle \pre(\dTwo), \pre(\dOne) \rangle}_B
        + \omega^2 \|\pre(\dTwo)\|^2
    \end{align*}
    Due to the bias introduced by nonlinear preconditioning, even after taking expectations in the equality above, we cannot immediately obtain ``good'' terms and so we will have to bound the terms $A$ and $B$ using different techniques. For $A$, we have the following:
    \begin{align*}
        A &= \underbrace{-2\omega \langle \pre(Fz^k),z^k-x^\star \rangle}_{\leq 0}-2\omega \langle \pre(\dTwo)-\pre(Fz^k),z^k-x^\star \rangle
        \\
        &\leq \omega^2\|z^k-x^\star\|^2 + \|\pre(\dTwo)-\pre(Fz^k)\|^2
        \\
        &\leq 2\omega^2 \|x^k-x^\star\|^2 + 2 \omega^2 \gamma^2 \|\pre(\dOne)\|^2 + \|\pre(\dTwo)-\pre(Fz^k)\|^2,
    \end{align*}
    where the nonnegativity of the inner product follows by \Cref{assump:generalized-wmvi} with $\rho=0$ and the first inequality follows by Young's inequality \eqref{eq:fenchel_young} with $e=1$. The second inequality follows by the $z^k$-update \eqref{eq:z_update} along with \eqref{eq:young} with $e=1$.

    Moving on to $B$ we get
    \begin{align*}
        B = \underbrace{\omega \gamma \|\pre(\dTwo) - \pre(\dOne)\|^2}_C - \omega \gamma\|\pre(\dTwo)\|^2 \underbrace{- \omega \gamma\|\pre(\dOne)\|^2}_D.
    \end{align*}
    The r.h.s.\ has one ``bad term'', $C$, that we will move on to bound using $\pre$-Lipschitz continuity of $F$ and two ``good'' terms from which we will obtain a negative term of the form $-\eta\|\pre(Fx^k)\|^2$ for some positive $\eta$, replacing the stochastic estimator with the quantity of interest. Adding and subtracting $\pre(Fz^k) - \pre(Fx^k)$ and applying \eqref{eq:young} twice with $e=1$ gives
    \begin{align*}
        C &\leq 2\omega \gamma \|\pre(Fz^k) - \pre(Fx^k)\|^2 + 4\omega \gamma\|\pre(\dTwo)-\pre(Fz^k)\|^2+4\omega \gamma\|\pre(Fx^k) - \pre(\dOne)\|^2
        \\
        & \leq 2\omega \gamma^3L^2\|\pre(\dOne)\|^2 + 4\omega \gamma\|\pre(\dTwo)-\pre(Fz^k)\|^2 +4\omega\gamma\|\pre(Fx^k) - \pre(\dOne)\|^2 ,
    \end{align*}
    where the second inequality follows by \cref{def:aniso_cont} between points $x^k$ and $z^k$ and substituting the $z$-update \eqref{eq:z_update}. For $D= -\tfrac{\omega\gamma}{2}\|\pre(\dOne)\|^2-\tfrac{\omega\gamma}{2}\|\pre(\dOne)\|^2 $, applying \eqref{eq:neg_young} to the first term, with $e=1$ and $a = \pre(Fx^k)$, $b = \pre(\dOne)-\pre(Fx^k)$ we get:
    \begin{align*}
        D \leq -\tfrac{\omega \gamma}{4} \|\pre(Fx^k)\|^2 + \tfrac{\omega \gamma}{2}\|\pre(\dOne)-\pre(Fx^k)\|^2 - \tfrac{\omega \gamma}{2}\|\pre(\dOne)\|^2.
    \end{align*}
    Putting everything together and denoting $e_k \coloneq \|\pre(\dOne)-\pre(Fx^k)\|^2$, $\lambda_k \coloneq \|\pre(\dTwo)-\pre(Fz^k)\|^2$, we obtain
    \begin{multline*}
        \|x^{k+1}-x^\star\|^2 \leq (1+2\omega^2)\|x^k-x^\star\|^2 + (2 \omega^2 \gamma^2 + 2\omega \gamma^3L^2 - \tfrac{\omega \gamma}{2})\|\pre(\dOne)\|^2
        \\
        + \omega(\omega-\gamma)\|\pre(\dTwo)\|^2 + \lambda_k + 4\omega \gamma \lambda_k + \tfrac{9}{2}\omega \gamma e_k -\tfrac{\omega \gamma}{4} \|\pre(Fx^k)\|^2.
    \end{multline*}
    By $\omega=\tfrac{\gamma}{\sqrt{K}}$ and $K\geq1$, we have $\omega\leq\gamma$, hence $\omega(\omega - \gamma) \leq 0$. Moreover, $2\omega\gamma + 2\gamma^2L^2 - \tfrac{1}{2} \leq 2\gamma^2 (1 + L^2) - \tfrac{1}{2} \leq 0$ due to $\gamma \leq \tfrac{1}{2\sqrt{1+L^2}}$. It follows that $2\omega^2\gamma^2 + 2\omega\gamma^3L^2 - \tfrac{\omega\gamma}{2} = \omega\gamma (2\omega\gamma + 2\gamma^2L^2 - \tfrac{1}{2}) \leq 0$. The above bound thus simplifies to
    \begin{align*}
        \|x^{k+1}-x^\star\|^2 \leq (1+2\omega^2)\|x^k-x^\star\|^2 + \lambda_k + 4\omega \gamma \lambda_k + \tfrac{9}{2}\omega \gamma e_k -\tfrac{\omega \gamma}{4} \|\pre(Fx^k)\|^2.
    \end{align*}
    We now have to bound $e_k, \lambda_k$, which we will do using the Lipschitz continuity of $\pre$. Due to the chosen batch sizes, and \eqref{eq:variance}, we have
    \[
        \bE[e_k \mid \cF_{k-1}] \leq \tfrac{L^2_\pre\sigma^2}{N_x} \leq \tfrac{L^2_\pre\sigma^2}{\sqrt{K}}, \qquad \bE[\lambda_k \mid \cG_k] \leq \tfrac{L^2_\pre\sigma^2}{N_z} = \tfrac{L^2_\pre\sigma^2}{K}.
    \]
    Taking expectations in the descent inequality and applying the tower property we are led to the following:
    \[
        \bE[\|x^{k+1}-x^\star\|^2] \leq (1+\tfrac{2\gamma^2}{K})\bE[\|x^k-x^\star\|^2] + L_{\pre}^2\tfrac{\sigma^2}{K} + \tfrac{9}{2}L_{\pre}^2 \tfrac{\gamma^2 \sigma^2}{K}
        + 4L_{\pre}^2\tfrac{\gamma^2 \sigma^2}{K^{3/2}} -\tfrac{\gamma^2}{4\sqrt{K}} \bE[\|\pre(Fx^k)\|^2].
    \]
    Using \cref{thm:sequences} with $V_k = \bE[\|x^k-x^\star\|^2]$, $\alpha = \tfrac{2\gamma^2}{K}$,  $u_k = \tfrac{\gamma^2}{4\sqrt{K}} \bE[\|\pre(Fx^k)\|^2]$ and $C = L_{\pre}^2\tfrac{\sigma^2}{K} + \tfrac{9}{2}L_{\pre}^2 \tfrac{\gamma^2 \sigma^2}{K}+ 4L_{\pre}^2\tfrac{\gamma^2 \sigma^2}{K^{3/2}}$ we get
    \[
        \frac{\gamma^2}{4\sqrt{K}} \sum_{k=0}^{K-1}\bE[\|\pre(Fx^k)\|^2] \leq \exp(2\gamma^2)\left[\|x^0-x^\star\|^2 + \frac{L_{\pre}^2}{2\gamma^2}\left( \frac{9}{2} \gamma^2 \sigma^2 + \sigma^2+ 4\frac{\gamma^2 \sigma^2}{K^{1/2}}\right)\right]
    \]
    and multiplying both sides by $\frac{4}{\gamma^2\sqrt{K}}$ yields the claimed result.
\end{proof}

\section{Proofs of \Cref{sec:OC}}
We first recall a standard result regarding the Fitzpatrick function of a maximally monotone operator, see for example \cite[p.\,761]{borwein2016maximal}.
\begin{proposition} \label{thm:fitz_prop}
    Let $B$ be maximally monotone. Then,
    \begin{equation} \label{eq:Fitz_lower_bound}
        f_B(x,y) \geq \langle y, x\rangle \qquad \forall x,y \in E.
    \end{equation}
    Moreover, $f_B(x,y) = \langle y, x\rangle$ if and only if $y \in Bx$.
\end{proposition}
\subsection{Proof of \cref{thm:var_form}}
\begin{proof}
    Since $F$ is $\nabla \phi^*$-Lipschitz continuous with constant $L$, the preconditioned operator $\nabla \phi^* \circ F$ is hypomonotone with constant $L$ and thus $B = \nabla\phi^* \circ F + L \id$ is maximally monotone by \cite[Corollary 20.28]{bauschke2017correction}. First, we show that for any $x \in X_{x_0}$, $J(x) \geq 0$. From \eqref{eq:Fitz_lower_bound} we have that
    \begin{align*}
        f_B(x, Lx-\dot x) - L\|x\|^2 \geq L\|x\|^2 - \langle x,\dot x \rangle - L\|x\|^2 = -\langle x,\dot x \rangle = -\tfrac{d}{dt}(\tfrac{1}{2}\|x\|^2).
    \end{align*}
    Integrating from $0$ to $T$ and using that $x$ is absolutely continuous we obtain
    \begin{equation*}
        \int_0^T [f_B(x, Lx-\dot x) - L\|x\|^2] dt \geq \tfrac{1}{2}\|x_0\|^2 - \tfrac{1}{2}\|x(T)\|^2
    \end{equation*}
    or, equivalently, $J(x) \geq 0$.

    The equality case of  \cref{thm:fitz_prop} shows that for $\bar x$ solving \eqref{eq:cont_time}, we have $J(\bar x) = 0$.

    Conversely, suppose $\bar x \in X_{x_0}$ is such that $J(\bar x) = 0$. Then,
    \begin{align*}
        \int_0^T [f_B(\bar x, L\bar x-\dot {\bar x}) - (L\|\bar x\|^2 - \langle \bar x, \dot {\bar x} \rangle)] dt = 0
    \end{align*}
    which, along with the fact that the integrand is nonnegative, implies that it vanishes a.e.\ on $(0,T)$. Hence, $f_B(\bar x, L\bar x-\dot {\bar x}) = L\|\bar x\|^2 - \langle \bar x, \dot {\bar x} \rangle = \langle L\bar x - \dot {\bar x},\bar x \rangle$ and from the equality case in \cref{thm:fitz_prop}, $ L\bar x-\dot {\bar x} = B\bar x$ or, equivalently, $\dot {\bar x}(t) = -\nabla \phi^*(F\bar x(t))$ almost everywhere on $(0,T)$.
\end{proof}

\subsection{Proof of \cref{thm:control}}
\begin{proof}
    To begin with, consider any stabilizing control $u$ and note that
    \begin{align*}
        \tfrac{d}{dt}V(x(t)) = \langle x(t)-x^\star,u(t) \rangle.
    \end{align*}
    Moreover, by \eqref{eq:Fitz_lower_bound}, we have for all $t \in [0,+\infty)$ that\[q(x(t),u(t)) = f_B(x(t),Lx(t)- u(t)) - L\|x(t)\|^2 + \langle u(t), x^\star\rangle \geq \langle x(t),-u(t) \rangle + \langle x^\star,u(t) \rangle,\] implying that $\tfrac{d}{dt}V(x(t)) + q(x(t),u(t)) \geq 0$ for almost every $t \geq 0$. We then have
    \begin{align*}
        V(x(t))-V(x_0) = \int_0^t \tfrac{d}{ds}V(x(s))ds \geq -\int_0^t q(x(s),u(s))ds.
    \end{align*}
    Now, since $u$ is stabilizing, $x(t)\to x^\star$ and hence $V(x(t))\to0$. Passing to the limit in the inequality above we get $J_\infty(x_0,u)=\int_0^\infty q(x(s),u(s))\,ds \geq V(x_0)$.

    Regarding the second statement, we first prove that $u$ is a stabilizing control. We have the following:
    \begin{align*}
        \tfrac{d}{dt}V(x(t)) = -\langle x(t)-x^\star, \pre(Fx(t))\rangle = -\isocoe(\|Fx(t)\|)\langle x(t)-x^\star,Fx(t) \rangle < 0,
    \end{align*}
    where the strict inequality holds whenever $x(t)\neq x^\star$ and follows from strict monotonicity of $F$, positivity of $\kappa$, and $Fx^\star = 0$. Moreover, $V$ is coercive, $\lim_{\|x\|\to \infty}V(x) = +\infty$, and smooth and strongly convex, and $V(x) = 0$ if and only if $x = x^\star $. Therefore, $V$ is a global control-Lyapunov function \cite[Definition 5.7.1]{sontag2013mathematical} and $x^\star$ is a globally asymptotically stable equilibrium, which in turn implies that $u$ is stabilizing. For this control, the equality characterization in \cref{thm:fitz_prop} gives $q(x(t),u(t))=-\tfrac{d}{dt}V(x(t))$. Integrating therefore yields $J_\infty(x_0,u) = V(x_0)$.
\end{proof}

\section{Experiments} \label{app:experiments}
Throughout this appendix, we consider the isotropic preconditioners
\[
    \pre(y)=y \textnormal{ (no preconditioning)}, \qquad \pre(y) =  \min\{1, \tfrac{\tau}{\|y\|}\}y \textnormal{ (Clip)}, 
\]
\[
    \pre(y) = y/(1+\|y\|/\tau) \textnormal{ (Softclip)}, \quad \pre(y) = \arsinh(\|y\|)y/\|y\| \textnormal{ (Hyper)}.
\]
At $0$, we define $\pre(0)=0$ for each of the preconditioners. The code for reproducing the experiments is publicly available\footnote{\url{https://github.com/JanQ/dual-space-prec-vis}}. The synthetic experiments and the WGAN experiment were run on CPU, while the MNIST GAN experiment was run on a single H100 GPU.

\subsection{Synthetic experiments}

We first evaluate the effect of nonlinear preconditioning on the \emph{PolarGame} and \emph{Forsaken} examples from \cite[Example 3 \& Example 6]{pethick2023escaping}. In both cases, we apply \eqref{eq:EG+} with $\beta=\tfrac12$ and compare no preconditioning (CEG+), with the Softclip (Softclip-CEG+) and Hyper (Hyper-CEG+) preconditioners.

\paragraph{Problem descriptions}
For PolarGame, we take $C=[-1.1,1.1]^2$ and
$F(x_1,x_2)=(q(x)x_1-x_2,q(x)x_2+x_1)$, where
$q(x)=(\|x\|^2-1)(16\|x\|^2-9)/16$.

For Forsaken, we take $C=[-1.5,1.5]^2$ and
$F(x_1,x_2)=(x_2-0.45+g(x_1),-x_1+g(x_2))$, where
$g(t)=\tfrac12t-2t^3+t^5$.

\paragraph{Tuning procedure}
For each problem, we sample $200$ initial points uniformly from $C$. We then search over $\gamma\in\{0.05,0.10,\ldots,2.00\}$ and, for Softclip-CEG+, additionally over $\tau\in\{1,2,5,10,20\}$. Each run of \eqref{eq:EG+} is allowed at most $2500$ iterations, and we stop when $\|x-P_C(x-Fx)\| \leq 10^{-6}$. We select the hyperparameters that maximize the success rate, breaking ties by minimizing first the median and then the $90$th percentile of the number of calls to $F$ used in the updates.

The selected parameters are subsequently evaluated on $2000$ independently sampled initial points, with a maximum of $5000$ iterations. The results are reported in \cref{tab:toy-experiments}. Both nonlinear preconditioners reduce the median number of calls to $F$ used in the updates.

\begin{table}[hbpt]
    \centering
    \setlength{\tabcolsep}{6pt}
    \renewcommand{\arraystretch}{1.05}
    \caption{Results on the PolarGame and Forsaken problems. Median calls are computed over successful runs only.}
    \label{tab:toy-experiments}
    \begin{tabular}{@{}llccrr@{}}
        \toprule
        Problem & Method & $\gamma$ & $\tau$ & Success & Median $F$ calls \\
        \midrule
        \multirow{3}{*}{PolarGame}
        & CEG+          & $0.95$ & --   & $99.30\%$  & $34$ \\
        & Softclip-CEG+ & $1.20$ & $10$ & $99.45\%$  & $\mathbf{22}$ \\
        & Hyper-CEG+    & $1.20$ & --   & $\mathbf{100.00\%}$ & $\mathbf{22}$ \\
        \midrule
        \multirow{3}{*}{Forsaken}
        & CEG+          & $1.10$ & --  & $99.25\%$ & $70$ \\
        & Softclip-CEG+ & $1.40$ & $5$ & $99.75\%$ & $\mathbf{54}$ \\
        & Hyper-CEG+    & $1.35$ & --  & $\mathbf{99.80\%}$ & $\mathbf{54}$ \\
        \bottomrule
    \end{tabular}
\end{table}

\subsection{Eight-Gaussian WGAN-GP experiment}
\label{subsec:eight-gaussian}
We next evaluate the effect of nonlinear preconditioning for training a Wasserstein GAN with gradient penalty (WGAN-GP) in the stochastic setting. Specifically, we aim to learn a mixture of Gaussians using \eqref{eq:EG+-variable_stoch}, and compare no preconditioning (SEG) with the Clip (Clip-SEG), Softclip (Softclip-SEG), and Hyper (Hyper-SEG) preconditioners.

\paragraph{Problem description}
Consider a mixture of eight two-dimensional Gaussians
\[
    p_{\mathrm{data}} = \frac18\sum_{j=0}^{7}\mathcal{N}(c_j,\sigma^2I_2), \qquad c_j=\bigl(r\cos(\pi j/4),r\sin(\pi j/4)\bigr)
\]
with $\sigma=0.08$ and $r=2$. We train a regularized Wasserstein GAN on this mixture, which is a commonly used benchmark for mode collapse \cite{mishchenko2020revisiting}. For the generator $G_{\phi}$ and critic $D_{\psi}$, we utilize MLPs with two hidden layers of width $32$ and $\tanh$ activations. Consider the two losses
\begin{align*}
    L_G(\phi, \psi) &= -\bE_{z\sim p_z}[D_\psi(G_\phi(z))] + \frac{\lambda}{2} \|\phi\|^2, \\
    L_D(\phi, \psi) &= \bE_{z\sim p_z} [D_\psi(G_\phi(z))] - \bE_{x\sim p_{\rm data}} [D_\psi(x)] + \lambda_{\rm gp} \bE_{x,z,\varepsilon}\left[(\|\nabla_{\hat{x}} D_\psi(\hat{x})\| - 1)^2 \right] + \frac{\lambda}{2} \|\psi\|^2
\end{align*}
where we take the latent distribution $p_z = \cN(0, I_2)$ and the regularization parameter $\lambda_{\rm gp} = 10^{-4}$ and $\lambda = 0.02$. In the third term of $L_D$, we have that $x\sim p_{\mathrm{data}}$, $z\sim p_z$, and $\varepsilon\sim\operatorname{Unif}[0,1]$ are sampled independently, and $\widehat{x}=\varepsilon x+(1-\varepsilon)G_\phi(z)$. The corresponding operator is
\[
    F(\phi, \psi) = \begin{pmatrix}
        \nabla_\phi L_G(\phi, \psi) \\ \nabla_\psi L_D(\phi, \psi)
    \end{pmatrix}.
\]

\paragraph{Tuning procedure}
We run each method for $T = 1024$ iterations with $\omega = \beta\gamma$, where $\beta=0.9$. At each iteration, we utilize $N_x = \sqrt{T} = 32$ and $N_z = T = 1024$. The stepsize $\gamma$ is initially $\gamma_0$ and is multiplied halfway through by a decay factor $\rho$. We perform a coarse search over $\gamma_0\in\{0.0125, 0.025, 0.05, 0.1, 0.2\}$ and $\rho\in\{0.2, 0.5, 1\}$. For Clip-SEG and Softclip-SEG, we additionally search over $\tau\in\{0.5, 1, 2\}$. Each coarse configuration is evaluated over ten random seeds using $2000$ generated samples, and parameter choices are ranked by their mean sliced Wasserstein distance (SWD). We retain the five best configurations and refine $\gamma$ and $\tau$, multiplying them by $\{2^{-1/2},1,2^{1/2}\}$. The refined parameters are then evaluated over ten seeds, and the choice with the smallest mean sliced Wasserstein distance is selected.

\paragraph{Evaluation}
The empirical sliced Wasserstein-1 distance is computed as follows: we sample $64$ random unit directions, then project both the generated and target samples onto each direction, and average the one-dimensional Wasserstein-1 distances over these projections.

To assess mode coverage, we additionally consider the distance of each generated point to its nearest Gaussian center. A generated point is classified as being near a mode if the distance is less than $0.35$, and we report the fraction of generated points satisfying this condition. We also report the number of covered modes, where a mode is considered covered if more than $1\%$ of all generated samples lie near this mode.

After tuning, we evaluate each algorithm on $20$ additional seeds, using $3000$ generated samples per run. The selected parameters and results are shown in \cref{tab:wgan_results}. Overall, the preconditioned variants improve distributional fit and reduce mode collapse, with the Hyper-SEG giving the strongest results across the three reported metrics.

\begin{table}[hbpt] 
    \centering
    \small
    \setlength{\tabcolsep}{3pt}
    \caption{Eight-Gaussian WGAN results: selected parameters and mean $\pm$ standard deviation for the SWD, covered modes, and near-mode fraction. SEG diverged for one seed, and the statistics are calculated over the $19$ remaining seeds instead.} \label{tab:wgan_results}
    \label{tab:eight-gaussian-results}
    \begin{tabular}{@{}lccc|ccc@{}}
        \toprule
        Method & $\gamma$ & Parameter $\tau$ & Decay $\rho$
        & SWD $\downarrow$ & Modes $\uparrow$ & Near-mode fraction $\uparrow$ \\
        \midrule
        SEG$^{\dagger}$ & $0.05\sqrt{2}$ & -- & $0.2$
        & $0.291 \pm 0.187$ & $6.05 \pm 2.07$ & $0.193 \pm 0.099$ \\
        Clip-SEG & $0.1\sqrt{2}$ & $2$ & $0.5$
        & $0.258 \pm 0.211$ & $6.65 \pm 2.37$ & $0.416 \pm 0.227$ \\
        Softclip-SEG & $0.1\sqrt{2}$ & $2\sqrt{2}$ & $0.5$
        & $0.205 \pm 0.073$ & $7.25 \pm 1.52$ & $0.332 \pm 0.098$ \\
        Hyper-SEG & $0.1\sqrt{2}$ & -- & $0.2$
        & $\mathbf{0.186} \pm 0.067$ & $\mathbf{7.55} \pm 0.83$
        & $\mathbf{0.420} \pm 0.147$ \\
        \bottomrule
    \end{tabular}
\end{table}

\subsection{MNIST GAN experiment}
Finally, we consider GAN training on the MNIST dataset \cite{deng2012mnist}. As in the preceding experiment, we compare stochastic extragradient without preconditioning (SEG) with the Clip (Clip-SEG), Softclip (Softclip-SEG), and Hyper (Hyper-SEG) preconditioners.

\paragraph{Problem description}
The generator $G_\phi$ is a convolutional network with $64$ base channels, group normalization and ReLU activations, followed by a $\tanh$ output layer. The discriminator $D_\psi$ consists of three convolutional layers with LeakyReLU activations followed by a linear output layer. We consider the loss functions
\begin{align*}
    L_G(\phi,\psi) &= \bE_{z\sim p_z}\left[\log\left(1+\exp\left(-D_\psi(G_\phi(z))\right)\right)\right], \\
    L_D(\phi,\psi) &= \bE_{x\sim p_{\mathrm{data}}} \left[\log\left(1+\exp\left(-D_\psi(x)\right)\right)\right] +\bE_{z\sim p_z} \left[\log\left(1+\exp\left(D_\psi(G_\phi(z))\right)\right)
    \right] \\
    &\qquad\qquad\qquad + \frac{\lambda_{\mathrm{out}}}{2} \left(\bE_{x\sim p_{\mathrm{data}}}[D_\psi(x)^2] + \bE_{z\sim p_z}[D_\psi(G_\phi(z))^2] \right) + \frac{\lambda_{\mathrm{wd}}}{2}\|\psi\|^2,
\end{align*}
where we take $p_z= \mathcal{N}(0, I_{128})$, $\lambda_{\rm out} = 10^{-3}$ and $\lambda_{\rm wd}=10^{-5}$. The corresponding operator is
\[
    F(\phi, \psi) = \begin{pmatrix}
        \nabla_\phi L_G(\phi, \psi) \\
        \nabla_\psi L_D(\phi, \psi)
    \end{pmatrix}.
\]

\paragraph{Tuning procedure}
As in the preceding experiment, we run each method for $T=1024$ iterations with $\omega=\beta\gamma$, where $\beta=0.9$. We also utilize minibatches of sizes $N_x=\sqrt{T}=32$ and $N_z = T =1024$. The stepsize $\gamma$ is initially $\gamma_0$ and multiplied halfway through by a decay factor $\rho$. We perform a coarse search over $\gamma_0\in \{0.00625, 0.0125, 0.025, 0.05, 0.1, 0.2\}$ and $\rho\in\{0.2, 0.5, 1\}$. For Clip-SEG and Softclip-SEG, we additionally search over $\tau\in\{0.5, 1, 2\}$. Each parameter setting is evaluated over five random seeds using $2000$ generated samples and the three settings with lowest MNIST Fréchet distance (mFID) are chosen. We then perform a refinement, multiplying them by $\{2^{-1/2}, 1, 2^{1/2}\}$. These are then evaluated over ten random seeds, again using $2000$ generated samples, and the parameters leading to the lowest mFID are chosen.

\paragraph{Evaluation}
The mFID is computed as follows: we train a separate convolutional classifier for MNIST digits, and extract a $128$-dimensional hidden representation for the real dataset, after which the mFID is computed as the Fréchet distance between fitted Gaussian distributions to the classifier features of the real and generated images. Additionally, we compute the number of covered digit classes using predictions from the same classifier. A class is considered covered if at least $1\%$ of all generated images are assigned to it.

\begin{figure}[hbpt]
    \centering

    \begin{subfigure}[t]{0.24\textwidth}
        \centering
        \includegraphics[width=\linewidth]{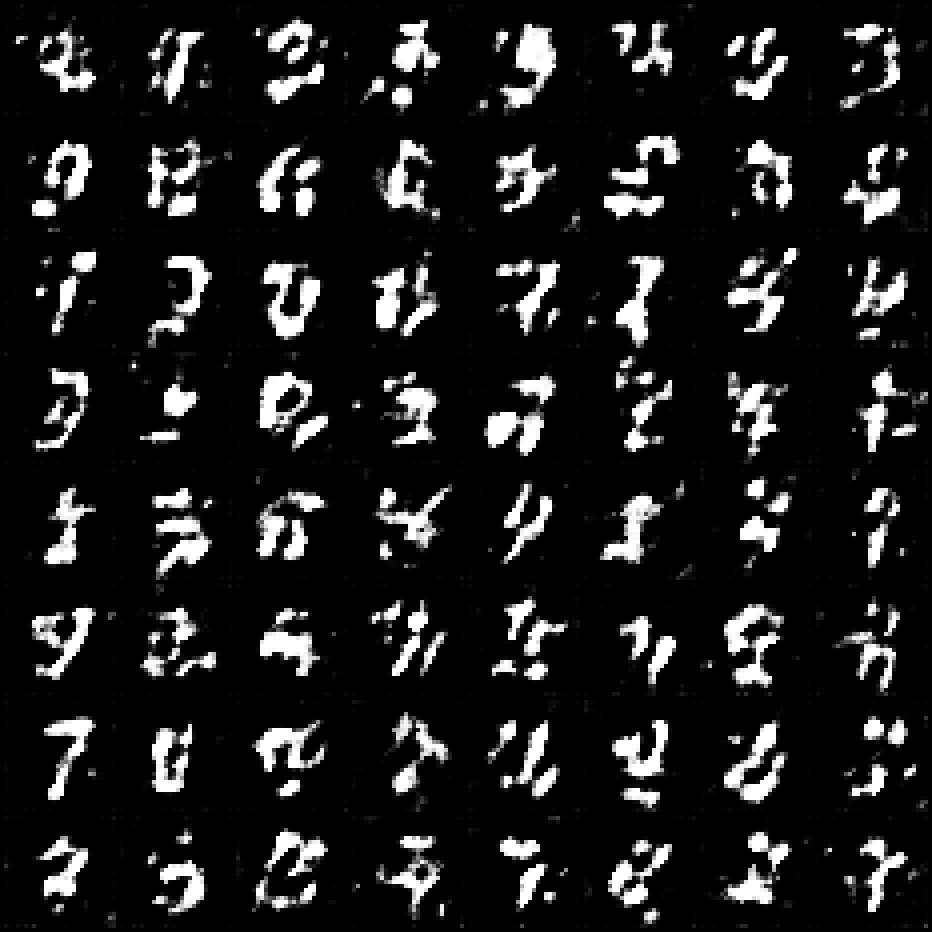}
        \caption{SEG}
    \end{subfigure}
    \hfill
    \begin{subfigure}[t]{0.24\textwidth}
        \centering
        \includegraphics[width=\linewidth]{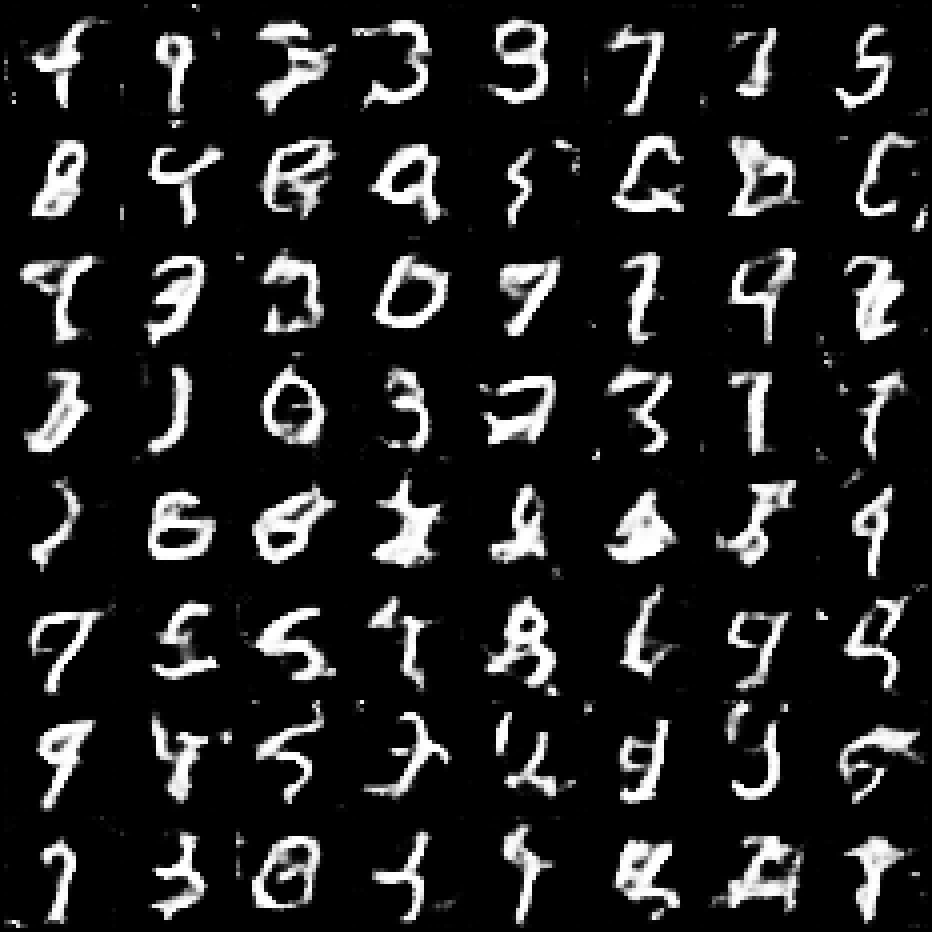}
        \caption{Clip-SEG}
    \end{subfigure}
    \hfill
    \begin{subfigure}[t]{0.24\textwidth}
        \centering
        \includegraphics[width=\linewidth]{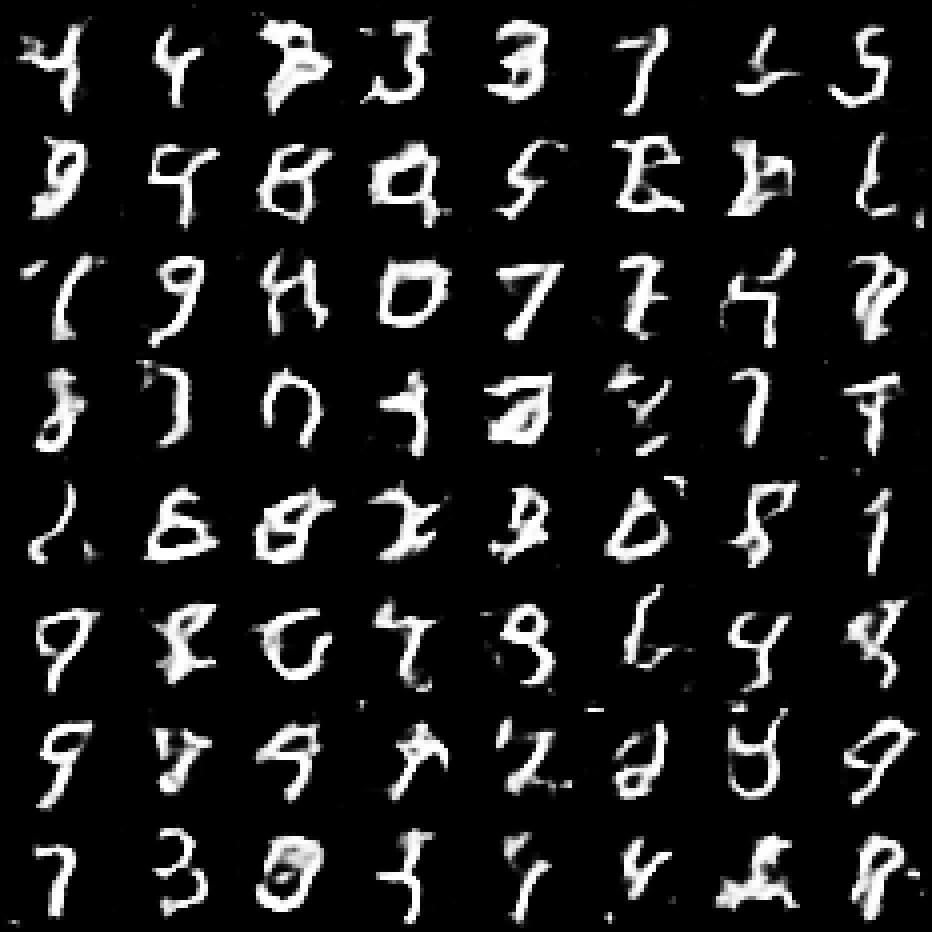}
        \caption{Softclip-SEG}
    \end{subfigure}
    \hfill
    \begin{subfigure}[t]{0.24\textwidth}
        \centering
        \includegraphics[width=\linewidth]{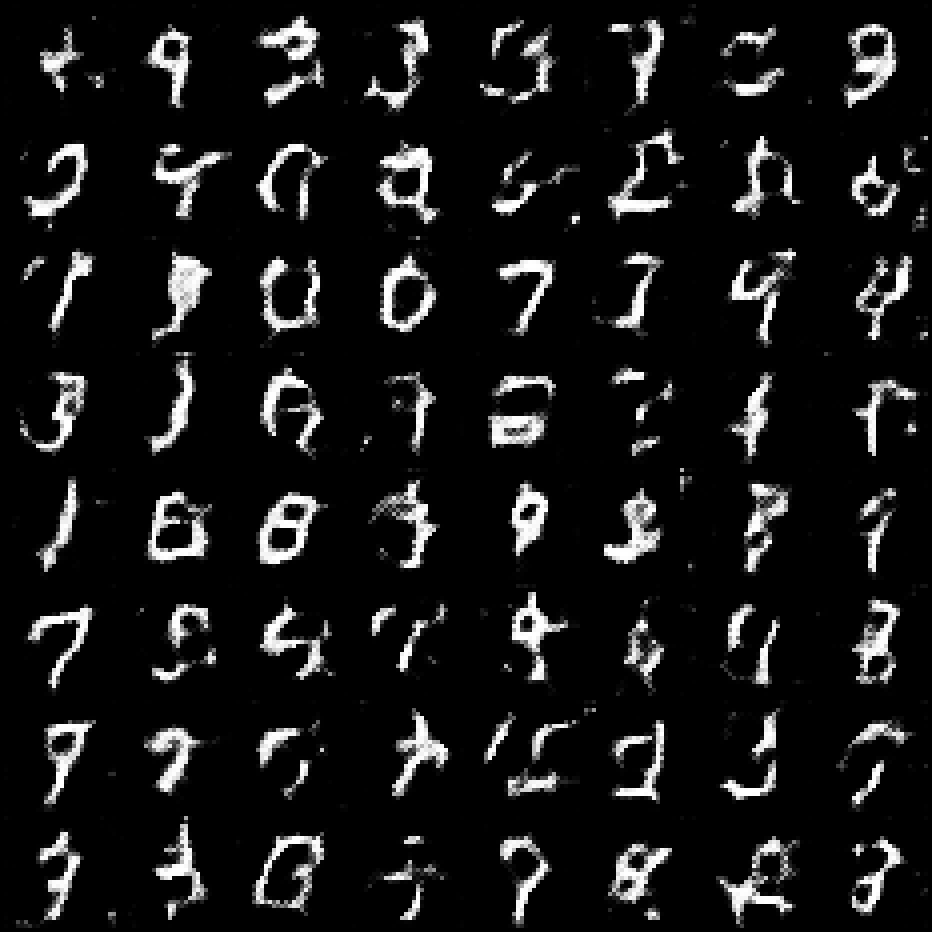}
        \caption{Hyper-SEG}
    \end{subfigure}

    \caption{Generated MNIST samples.}
    \label{fig:mnist-generated-samples}
\end{figure}

After tuning, we evaluate each algorithm on five additional seeds, using $5000$ generated samples per run. The selected parameters and results are shown in \cref{tab:mnist_gan_results}, and some representative generated images are shown in \cref{fig:mnist-generated-samples}. The results indicate that the preconditioned variants consistently outperform the nonpreconditioned SEG baseline, both in the mFID and the quality of the generated images.

\begin{table}[hbpt]
    \centering
    \small
    \setlength{\tabcolsep}{6pt}
    \caption{MNIST GAN results: selected parameters and mean $\pm$
    sample standard deviation of mFID and class coverage over five runs.}
    \label{tab:mnist_gan_results}
    \begin{tabular}{@{}lccc|cc@{}}
        \toprule
        Method & $\gamma$ & Parameter $\tau$ & Decay $\rho$
        & mFID $\downarrow$ & Class coverage $\uparrow$ \\
        \midrule
        SEG & $0.0125$ & -- & $0.2$
        & $485.89 \pm 37.08$
        & $9.80 \pm 0.45$ \\
        Clip-SEG & $0.1\sqrt{2}$ & $2$ & $1$
        & $161.41 \pm 30.31$
        & $\mathbf{10.00} \pm 0.00$ \\
        Softclip-SEG & $0.2\sqrt{2}$ & $\sqrt{2}$ & $1$
        & $\mathbf{158.63} \pm 49.79$
        & $\mathbf{10.00} \pm 0.00$ \\
        Hyper-SEG & $0.05\sqrt{2}$ & -- & $1$
        & $306.66 \pm 25.77$
        & $\mathbf{10.00} \pm 0.00$ \\
        \bottomrule
    \end{tabular}
\end{table}

\end{document}